\documentclass{article}
\usepackage{amsmath}
\usepackage{amssymb}
\usepackage{amsthm}
\usepackage{pxfonts}
\usepackage{times}
\usepackage{bbm}
\usepackage{graphicx}
\usepackage{algorithm}
\usepackage{algorithmic}
\usepackage{clrscode}
\usepackage{fullpage}

\newcounter{thmnum}

\newtheorem{cor}[thmnum]{Corollary}
\newtheorem{thm}[thmnum]{Theorem}
\newtheorem{lem}[thmnum]{Lemma}

\newcommand{\polylog}{\mathop{\rm polylog}}

\begin{document}

% Title.
% ------
\title{On the Design of Deterministic Matrices for Fast Recovery of Fourier Compressible Functions}

\author{J. Bailey \thanks{Mathematics Department, Kansas State University, Manhattan, KS 66506, {\tt jamespb@ksu.edu}.%, research supported in part by a McNair Scholarship and an I-Center Undergraduate Scholarship from Kansas State University.
}
 \and\ 
M. A. Iwen\thanks{Contact Author.  Mathematics Department, Duke University, Durham, NC 27708, {\tt markiwen@math.duke.edu}, research supported in part by ONR N00014-07-1-0625 and NSF DMS-0847388.} \and\ 
C. V. Spencer\thanks{Mathematics Department, Kansas State University, Manhattan, KS 66506, {\tt cvs@math.ksu.edu}, research supported in part by NSA Young Investigators Grant H98230-10-1-0155.}} 

\maketitle
\thispagestyle{empty}

\begin{abstract}
We present a general class of compressed sensing matrices which are then demonstrated to have associated sublinear-time sparse approximation algorithms.
We then develop methods for constructing specialized matrices from this class which are sparse when multiplied with a discrete Fourier transform matrix.  Ultimately, these considerations improve previous sampling requirements for deterministic sparse Fourier transform methods.  
\end{abstract}

\section{Introduction}
\label{sec:intro}

This paper considers methods for designing matrices which yield near-optimal nonlinear approximations to the Fourier transform of a given function, $f: [0, 2\pi] \rightarrow \mathbbm{C}$.  Suppose that $f$ is a bandlimited function so that $\hat{f} \in \mathbbm{C}^{N}$, where $N$ is large.  An optimal $k$-term trigonometric approximation to $f$ is given by
$$f^{\rm opt}_k (x) = \sum^k_{j = 1} \hat{f} \left( \omega_j \right) \mathbbm{e}^{\mathbbm{i} \omega_j x},$$
where $\omega_1, \dots, \omega_{N} \in (-N/2, N/2] \cap \mathbbm{Z}$ are ordered by the magnitudes of their Fourier coefficients so that
$$\big| \hat{f}(\omega_1) \big| \geq \big| \hat{f}( \omega_2 ) \big| \geq \dots \geq \big| \hat{f}(\omega_{N}) \big|.$$
The optimal $k$-term approximation error is then
\begin{equation}
\left\| f - f^{\rm opt}_k \right\|_2 = \left\| \hat{f} - \hat{f}^{~\rm opt}_k \right\|_2.
\label{eqn:Aerror}
\end{equation}
It has been demonstrated recently that any periodic function, $f: [0, 2\pi] \rightarrow \mathbbm{C}$, can be accurately approximated via \textit{sparse Fourier transform (SFT)} methods which run in $O(k^2 \log^4 N)$ time (see \cite{JournalFFT, ImpFourier} for details).  When the function is sufficiently Fourier compressible (i.e., when $k << N$ yields a small approximation error in Equation~(\ref{eqn:Aerror}) above), these methods can accurately approximate $f$ much more quickly than traditional Fast Fourier Transform (FFT) methods which run in $O(N \log N)$ time.  Furthermore, these SFT methods require only $O(k^2 \log^4 N)$ function evaluations as opposed to the $N$ function evaluations required by a standard FFT method.

Although the the theoretical guarantees of SFT algorithms appear promising, current algorithmic formulations suffer from several practical shortfalls.  Principally, the algorithms currently utilize number theoretic sampling sets which are constructed in a suboptimal fashion.  In this paper we address this deficiency by developing computational methods for constructing number theoretic matrices of the type required by these SFT methods which are nearly optimal in size.  In the process, we demonstrate that this specific problem is a more constrained instance of a much more general matrix design problem with connections to compressed sensing matrix constructions \cite{Coherence, CS2, NearOpt, CS4, RIP1expand}, discrete uncertainty principles \cite{DonStark}, nonadaptive group testing procedures \cite{FGroupTest, MyGroupTest}, and codebook design problems \cite{CodebookLowerB, CodebookDesign1, CodebookDesign2} in signal processing.

\subsection{General Problem Formulation:  Compressed Sensing in the Fourier Setting}
\label{sec:GenProb}

Over the past several years, a stream of work in compressed sensing has provided a general theoretical framework for approximating general functions in terms of their optimal $k$-term approximation errors (see \cite{HolgerCSBook} and references therein).  Indeed, the SFT design problem we are considering herein also naturally falls into this setting.  Consider the following discretized version of the sparse Fourier approximation problem above:  Let $\vec{f} \in \mathbbm{C}^N$ be a vector of $N$ equally spaced evaluations of $f$ on $[0, 2\pi]$, and define $\mathcal{F}$ to be the $N \times N$ \textit{Discrete Fourier Transform (DFT)} matrix defined by $\mathcal{F}_{i,j} = \frac{\mathbbm{e}^{\frac{-2 \pi \mathbbm{i} \cdot i \cdot j}{N}}}{\sqrt{N}}$.  Note that $\mathcal{F} \vec{f}$ will be compressible (i.e., sparse).  Compressed sensing methods allow us to construct an $m \times N$ matrix, $\mathcal{M}$, with $m$ minimized as much as possible subject to the constraint that an associated approximation algorithm, $\Delta_{\mathcal{M}}: \mathbbm{C}^m \rightarrow \mathbbm{C}^N$, can still accurately approximate any given $\hat{f} = \mathcal{F} \vec{f}$ (and, therefore, $f$ itself).  More exactly, compressed sensing methods allow us to minimize $m$, the number of rows in $\mathcal{M}$, as a function of $k$ and $N$ such that
\begin{equation}
\left\| ~ \Delta_{\mathcal{M}} \left( \mathcal{M} \hat{f} \right) - \hat{f} ~ \right\|_p \leq C_{p,q} \cdot k^{\frac{1}{p} - \frac{1}{q}} \left\| ~ \hat{f} - \hat{f}^{~\rm opt}_k ~ \right\|_q
\label{eqn:Aerror2}
\end{equation}
holds for all $\hat{f} \in \mathbbm{C}^N$ in various fixed $l^p$,$l^q$ norms, $1 \leq q \leq p \leq 2$, for an absolute constant $C_{p,q} \in \mathbbm{R}$ (e.g., see \cite{BestkTerm, HolgerCSBook}).  Note that this implies that $\hat{f}$ will be recovered exactly if it contains only $k$ nonzero Fourier coefficients.  Similarly, it will be accurately approximated by $\Delta_{\mathcal{M}} \left( \mathcal{M} \hat{f} \right)$ any time it is well represented by its largest $k$ Fourier modes.

In this paper we will focus on constructing $m \times N$ compressed sensing matrices, $\mathcal{M}$, for the Fourier recovery problem which meet the following four design requirements:
\begin{enumerate}
\item \textbf{Small Sampling Requirements:} $\mathcal{M}\mathcal{F}$ should be highly column-sparse (i.e., the number of columns of $\mathcal{M}\mathcal{F}$ which contain nonzero entries should be significantly smaller than $N$).  Note that whenever $\mathcal{M}\mathcal{F}$ has this property we can compute $\mathcal{M} \hat{f}$ by reading only a small fraction of the entries in $\vec{f}$.  Once the number of required function samples/evaluations is on the order of $N$, a simple fast Fourier transform based approach will be difficult to beat computationally.
\item \textbf{Accurate Approximation Algorithms:} The matrix $\mathcal{M}$ needs to have an associated approximation algorithm, $\Delta_{\mathcal{M}}$, which allows accurate recovery.  More specifically, we will require an instance optimal error guarantee along the lines of Equation~(\ref{eqn:Aerror2}).
\item \textbf{Efficient Approximation Algorithms:} The matrix $\mathcal{M}$ needs to have an associated approximation algorithm, $\Delta_{\mathcal{M}}$, which is computationally efficient.  In particular,
the algorithm should be at least polynomial time in $N$ (preferably, $o(N \log N)$-time since $N$ is presumed to be large and we have the goal in mind of competing with an FFT).
\item \textbf{Guaranteed Uniformity:}  Given only $k, N \in \mathbbm{Z}^{+}$ and $p, q \in [1,2]$, one fixed matrix $\mathcal{M}$ together with a fixed approximation algorithm $\Delta_{\mathcal{M}}$ should be guaranteed to satisfy the three proceeding properties uniformly for all vectors $\hat{f} \in \mathbbm{C}^N$.
\end{enumerate}

The remainder of this paper is organized as follows:  We begin with a brief survey of recent sparse Fourier approximation techniques related to compressed sensing in Section~\ref{sec:RipSurvey}.  In Section~\ref{sec:BinaryIncoherent} we introduce matrices of a special class which are useful for fast sparse Fourier approximation and investigate their properties.  Most importantly, we demonstrate that any matrix from this class can be used in combination with an associated fast approximation algorithm in order to produce a sublinear-time (in $N$) compressed sensing method.  Next, in Section~\ref{sec:makeCoherent}, we present a deterministic construction of these matrices that specifically supports sublinear-time Fourier approximation.  In Section~\ref{sec:OptBuild} this matrix construction method is cast as an optimal design problem whose objective is to minimize Fourier sampling requirements.  Furthermore, lemmas are proven which allow the optimal design problem to be subsequently formulated as a linear integer program in Section~\ref{sec:LinIntProg}.  Finally, in Section~\ref{sec:Experiments}, we empirically investigate the sizes of the optimized deterministic matrices presented herein.

\section{Compressed Sensing and The Restricted Isometry Property}
\label{sec:RipSurvey}

Over the past few years, compressed sensing has focused primarily on utilizing matrices, $\mathcal{M}$, which satisfy the \textit{Restricted Isometry Principle (RIP)} \cite{CS2} in combination with $l^1$-minimization based approximation methods \cite{CS2, NearOpt, CS4}.  In fact, RIP matrices appear to be the critical partner in the RIP matrix/$l^1$-minimization pair since RIP matrices can also be used for compressed sensing with numerous other approximation algorithms besides $l^1$-minimization (e.g., Regularized Orthogonal Matching Pursuit \cite{ROMP,ROMPstable}, CoSaMP \cite{COSAMP}, Iterative Hard Thresholding \cite{HardThreshforCS}, etc.).  Hence, we will consider RIP matrices in isolation.

\newtheorem{Definition}{Definition}
\begin{Definition}
Let $p \in [1,\infty)$, $N,k \in \mathbbm{N}$, and $\epsilon \in (0,1)$.  A matrix $\mathcal{M}$ with complex entries has the Restricted Isometry Property, RIP$_{p}$($N$,$k$,$\epsilon$), if
$$ \left(1-\epsilon \right) \| x \|^p_{p} \leq \| \mathcal{M} x \|^p_{p} \leq \left(1+\epsilon\right) \| x \|^p_p$$
for all $x \in \mathbbm{C}^{N}$ containing at most $k$ nonzero coordinates.
\label{def:RIP}
\end{Definition}

It has been demonstrated that Fourier RIP$_{2}$($N$,$k$,$\epsilon$) matrices of size $O\left(k \frac{\log^4 N}{\epsilon^2} \right) \times N$ exist \cite{CSb1}.  More specifically, an $m \times N$ submatrix of the $N \times N$ \textit{Inverse DFT (IDFT)} matrix, $\mathcal{F}^{-1}$, formed by randomly selecting $m$ rows of $\mathcal{F}^{-1}$ will satisfy the RIP$_{2}$($N$,$k$,$\epsilon$) with high probability whenever $m$ is $\Omega\left(k \log^2 N \frac{\log^2 k}{\epsilon^2} \right)$ \cite{HolgerStableRecov}.  Such a matrix will clearly satisfy our small sampling requirement since any $m \times N$ submatrix of the $N \times N$ IDFT matrix will generate a vector containing exactly $m$ ones after being multiplied against the $N \times N$ DFT matrix.  Furthermore, $l^1$-minimization will yield accurate approximation of Fourier compressible signals when utilized in conjunction with an IDFT submatrix that has the RIP$_2$.  However, these random Fourier RIP$_2$ constructions have two deficiencies:  First, all existing approximation algorithms, $\Delta_{\mathcal{M}}$, associated with Fourier RIP$_{2}$($N$,$k$,$\epsilon$) matrices, $\mathcal{M}$, run in $\Omega \left( N \log N \right)$ time.  Thus, they cannot generally compete with an FFT computationally.  Second, randomly generated Fourier submatrices are only guaranteed to have the RIP$_2$ with high probability, and there is no tractable means of verifying that a given matrix has the RIP$_2$.  In order to verify Definition~\ref{def:RIP} for a given $m \times N$ matrix one generally has to compute the condition numbers of all ${N \choose k}$ of its $m \times k$ submatrices.

Several deterministic RIP$_{2}$($N$,$k$,$\epsilon$) matrix constructions exist which simultaneously address the guaranteed uniformity requirement while also guaranteeing small Fourier sampling needs \cite{DetFRIP, BestDRIP}.  However, they all utilize the notion of \textit{coherence} \cite{Coherence} which is discussed in Section~\ref{sec:Coherence}.  Hence, we will postpone a more detailed discussion of these methods until later.  For now, we simply note that no existing deterministic RIP$_2$($N$,$k$,$\epsilon$) matrix constructions currently achieve a number of rows (or sampling requirements), $m$, that are $o \left(k^2 \polylog(N) \right)$ for all $k = o\left(\sqrt{N}\right)$ as $N$ grows large.  In contrast, RIP matrix constructions related to highly unbalanced expander graphs can currently break this ``quadratic-in-$k$ bottleneck''.

\subsection{Unbalanced Expander Graphs}

Recently it has been demonstrated that the rescaled adjacency matrix of any unbalanced expander graph will be a RIP$_{1}$ matrix \cite{RIP1expand, RIP2Expand2}.  

\begin{Definition}
Let $N,k,d \in \mathbbm{N}$, and $\epsilon \in (0,1)$.  A simple bipartite graph $G = (A,B,E)$ with $|A| \geq |B|$ and left degree at least $d$ is a $(k,d,\epsilon)$-unbalanced expander if, for any $X \subset A$ with $|X| \leq k$, the set of neighbors, $|N(X)|$, of $X$ has size $|N(X)| \geq (1-\epsilon)d|X|$.
\label{def:Expander}
\end{Definition}

\newtheorem{Theorem}{Theorem}
\begin{Theorem} \textrm{\textbf{(See \cite{RIP1expand, RIP2Expand2}).}}
Consider an $m \times N$ matrix $\mathcal{M}$ that is the adjacency matrix of a regular $(k,d,\epsilon)$-unbalanced expander, where $1/\epsilon$ and $d$ are both smaller than $N$. Then, there exists an absolute constant $C > 1$ such that the rescaled matrix, $\mathcal{M}$ / $d^{1/p}$, satisfies the RIP$_{p}$($N$,$k$,$C\epsilon$) for all $1 \leq p \leq 1 ~+~1 / \log N$.
\label{thm:RIPexpand}
\end{Theorem}

Note that the RIP$_1$ property for unbalanced expanders is with respect to the $\textit{l}^1$ norm, not the $\textit{l}^2$ norm.  Nevertheless, matrices with the RIP$_{1}$ property also have associated approximation algorithms that can produce accurate sparse approximations along the lines of Equation~(\ref{eqn:Aerror2}).  Examples include $l^1$-minimization \cite{RIP1expand, RIP2Expand2} and Matching Pursuit \cite{IndykMP}.  Perhaps most impressive among the approximation algorithms associated with unbalanced expander graphs are those which appear to run in $o(N \log N)$-time (see the appendix of \cite{RIP2Expand2}).  Considering these results with respect to the four design requirements from Section~\ref{sec:GenProb}, we can see that expander based RIP methods are poised to satisfy both the second and third requirements.  Furthermore, by combining Theorem~\ref{thm:RIPexpand} with recent explicit constructions of unbalanced expander graphs \cite{Umans}, we can obtain an explicit RIP$_1$ matrix construction of near-optimal dimensions (which, among other things, shows that RIP$_1$ matrices may also satisfy our fourth Section~\ref{sec:GenProb} design requirement regarding guaranteed uniformity).  We have the following theorem:

\begin{Theorem}
Let $\epsilon \in (0,1)$, $p \in [1,1+1 / \log N]$, and $N,k \in \mathbbm{N}$ such that $N$ greater than both $1/\epsilon$ and $k$.  Next, choose any constant parameter $\alpha \in \mathbbm{R}^+$.  Then, there exists a constant $c \in \mathbbm{R}^+$ such that a $$O\left( k^{1+\alpha} \left( \log N \log k / \epsilon \right)^{2 + 2 / \alpha} \right) \times N$$ 
matrix guaranteed to have the RIP$_{p}$($N$,$k$,$\epsilon$) can be constructed in $O \left(N \cdot \left(\log N / \epsilon\right)^{c(1+1/\alpha)} \right)$-time.
\label{thm:ExplicitRIP1}
\end{Theorem}

\noindent \textit{Proof:}  Consider Theorem 1.3 in \cite{Umans} in combination with Theorem~\ref{thm:RIPexpand} above.~~$\Box$ \\

Theorem~\ref{thm:ExplicitRIP1} demonstrates the existence of deterministically constructible RIP$_1$ matrices with a number of rows, $m$, that scales like $O \left(k^{1 + \alpha} \polylog(N) \right)$ for all $k < N$ and fixed $\epsilon, \alpha \in (0,1)$.  Furthermore, the run time complexity of the RIP$_1$ construction algorithm is modest (i.e., $O\left(N^2 \right)$-time).  Although a highly attractive result, there is no guarantee that Guruswami et al.'s unbalanced expander graphs will generally have adjacency matrices, $\mathcal{M}$, which are highly column-sparse after multiplication against a DFT matrix (see design requirement number 1 in Section~\ref{sec:GenProb}).\footnote{In fact, $\mathcal{M}$ multiplied against a DFT matrix need not be exactly sparse.  By appealing to ideas from \cite{GeneralDeviants}, one can see that it is enough to have a relatively small perturbation of $\mathcal{M}$ be column-sparse after multiplication against a DFT matrix.}  Hence, it is unclear whether expander graph based RIP$_1$ results can be utilized to make progress on our compressed sensing matrix design problem in the Fourier setting.  Nevertheless, this challenging avenue of research appears potentially promising, if not intractably difficult.

\subsection{Incoherent Matrices}
\label{sec:Coherence}
As previously mentioned, all deterministic RIP$_{2}$($N$,$k$,$\epsilon$) matrix constructions (e.g., see \cite{DetCS1,DetFRIP,HolgerStableRecov,BestDRIP} and references therein) currently utilize the notion of coherence \cite{Coherence}.

\begin{Definition}
Let $\mu \in [0,1]$.  An $m \times N$ matrix, $\mathcal{M}$, with complex entries is called $\mu$-coherent if both of the following properties hold:
\begin{enumerate}
\item Every column of $\mathcal{M}$, denoted $\mathcal{M}_{\cdot,j} \in \mathbbm{C}^m$ for $0 \leq j \leq N-1$, is normalized so that $\| \mathcal{M}_{\cdot,j} \|_2 = 1$.
\item For all $j, l \in [0,N)$ with $j \neq l$, the associated columns $\mathcal{M}_{\cdot,j},\mathcal{M}_{\cdot,l} \in \mathbbm{C}^m$ have $\left| \mathcal{M}_{\cdot,j} \cdot \mathcal{M}_{\cdot,l} \right| \leq \mu$.
\end{enumerate}
\label{def:Coherence}
\end{Definition}

\begin{Theorem} \textrm{\textbf{(See \cite{HolgerStableRecov}).}}
Suppose that an $m \times N$ matrix, $\mathcal{M} \in \mathbbm{C}^{m \times N}$, is $\mu$-coherent.  Then, $\mathcal{M}$ will also have the RIP$_{2}$($N$,$k$,$(k-1) \mu$).
\label{thm:RIPcoherent}
\end{Theorem}

Matrices with small coherence are of interest in numerous coding theoretic settings.  Note that the column vectors of a real valued matrix with small coherence, $\mu$, collectively form a spherical code.  More generally, the columns of an incoherent complex valued matrix can be used to form codebooks for various channel coding applications in signal processing \cite{CohApp1,Grassman}.  These applications have helped to motivate a considerable amount of work with incoherent codes (i.e., incoherent matrices) over the past several decades.  As a result, a plethora of $\mu$-coherent matrix constructions exist (e.g., see \cite{CodebookLowerB,CodebookDesign1,CodebookDesign2,BestDRIP}, and references therein).

As we begin to demonstrate in the next section, matrices with low coherence can satisfy all four Fourier design requirements listed in Section~\ref{sec:GenProb}.  However, there are trade-offs.  Most notably, the Welch bound \cite{CodebookLowerB} implies that any $\mu$-coherent $m \times N$ matrix, $\mathcal{M} \in \mathbbm{C}^{m \times N}$, must have a number of rows
$$m \geq \frac{N}{(N-1) \mu^2 + 1}.$$
As a consequence, arguments along the lines of Theorem~\ref{thm:RIPcoherent} can only use $\mu$-coherent matrices to produce RIP$_{2}$($N$,$k$,$\epsilon$) matrices having $m = \Omega \left( k^2 / \epsilon^2 \right)$ rows. In contrast, $O\left(k \frac{\log^4 N}{\epsilon^2} \right) \times N$ Fourier RIP$_{2}$($N$,$k$,$\epsilon$) matrices are known to exist (see above).  Hence, although $\mu$-coherent matrices do allow one to obtain small Fourier sampling requirements, these sampling requirements all currently scale \textit{quadratically} with $k$ instead of \textit{linearly}.\footnote{It is worth noting that Bourgain et al. recently used methods from additive combinatorics in combination with modified coherence arguments to construct explicit $m \times N$ matrices, with $m = O(k^{2 - \epsilon'})$, which have the Fourier RIP$_{2}$($N$,$k$,$m^{-\epsilon'}$) whenever $k = \Omega(N^{1/2 - \epsilon'})$ \cite{BestDRIP}.  Here $\epsilon' > 0$ is some constant real number.  Hence, it is possible to break the previously mentioned ``quadratic bottleneck'' for RIP$_{2}$($N$,$k$,$\epsilon$) matrices when $k$ is sufficiently large.}

Setting aside the quadratic scaling of $m$ with $k$, we can see that several existing deterministic RIP$_{2}$($N$,$k$,$\epsilon$) matrix constructions based on coherence arguments (e.g., \cite{DetFRIP,BestDRIP}) immediately satisfy all but one of the Fourier design requirements listed in Section~\ref{sec:GenProb}.  First, these constructions lead to Fourier sampling requirements which, although generally quadratic in the sparsity parameter $k$, are nonetheless $o(N)$.  Second, these matrices can be used in conjunction with accurate approximation algorithms (e.g., $l^1$-minimization) since they will have the RIP$_{2}$.  Third, the deterministic nature of these RIP$_{2}$ matrices guarantees uniform approximation results for all possible periodic functions.  The only unsatisfied design requirement pertains to the computational efficiency of the approximation algorithms (see requirement 3 in Section~\ref{sec:GenProb}).  As mentioned previously, all existing approximation algorithms associated with Fourier RIP$_{2}$($N$,$k$,$\epsilon$) matrices run in $\Omega \left( N \log N \right)$ time.  In the next section we will present a general class of incoherent matrices which have fast approximation algorithms associated with them.  As a result, we will develop a general framework for constructing fast sparse Fourier algorithms which are capable of approximating compressible signals more quickly than standard FFT algorithms.

\section{A Special Class of Incoherent Matrices}
\label{sec:BinaryIncoherent}

In this section, we will consider binary incoherent matrices, $\mathcal{M} \in \{0,1\}^{m \times N}$, as a special subclass of incoherent matrices.  As we shall see, binary incoherent matrices can be used to construct RIP$_{2}$ matrices (e.g., via Theorem~\ref{thm:RIPcoherent}), unbalanced expander graphs (and, therefore, RIP$_{p \approx 1}$ matrices via Theorem~\ref{thm:RIPexpand}), and nonadaptive group testing matrices \cite{FGroupTest}.  In addition, we prove that any binary incoherent matrix can be modified to have an associated accurate approximation algorithm, $\Delta_{\mathcal{M}}: \mathbbm{C}^m \rightarrow \mathbbm{C}^N$, with sublinear $o(N)$ run time complexity.  This result generalizes the fast sparse Fourier transforms previously developed in \cite{ImpFourier} to the standard compressed sensing setup while simultaneously providing a framework for the subsequent development of similar Fourier results.  We will begin this process by formally defining $\left( K,\alpha \right)$-coherent matrices and then noting some
accompanying bounds. 

\begin{Definition}
Let $K,\alpha \in [1,m] \cap \mathbbm{N}$.  An $m \times N$ binary matrix, $\mathcal{M} \in \{0,1\}^{m \times N}$, is called $\left( K,\alpha \right)$-coherent if both of the following properties hold:
\begin{enumerate}
\item Every column of $\mathcal{M}$ contains at least $K$ nonzero entries.
\item For all $j, l \in [0,N)$ with $j \neq l$, the associated columns, $\mathcal{M}_{\cdot,j} ~and~ \mathcal{M}_{\cdot,l} \in \{0,1\}^m$, have $\mathcal{M}_{\cdot,j} \cdot \mathcal{M}_{\cdot,l} \leq \alpha$.
\end{enumerate}
\label{def:BinaryCoherence}
\end{Definition}

Several deterministic constructions for $\left( K,\alpha \right)$-coherent matrices have been implicitly developed as part of RIP$_2$ matrix constructions (e.g., see \cite{DetCS1,DetFRIP}).  It is not difficult to see that any $\left( K,\alpha \right)$-coherent matrix will be $\frac{\alpha}{K}$-coherent after having its columns normalized.  Hence, the Welch bound also applies to $\left( K,\alpha \right)$-coherent matrices.  Below we will both develop tighter lower row bounds, and provide a preliminary demonstration of the existence of fast $o(N)$-time compressed sensing algorithms related to incoherent matrices.  This will be done by demonstrating the relationship between $\left( K,\alpha \right)$-coherent matrices and group testing matrices.

\subsection{Group Testing: Lower Bounds and Fast Recovery}
\label{sec:GTcohMat}

Group testing generally involves the creation of testing procedures which are designed to identify a small number of interesting items hidden within a much larger set of uninteresting items \cite{FGroupTest,MyGroupTest}.  Suppose we are given a collection of $N$ items, each of which is either interesting or uninteresting.  The status of each item in the set can then be represented by a boolean vector $\vec{x} \in \{ 0,1 \}^N$.  Interesting items are denoted with a $1$ in the vector, while uninteresting items are marked with a $0$.  Because most items are uninteresting, $\vec{x}$ will contain at most a small number, $d < N$, of ones.  Our goal is to correctly identify the nonzero entries of $\vec{x}$, thereby recovering $\vec{x}$ itself.

Consider the following example.  Suppose that $\vec{x}$ corresponds to a list of professional athletes, at most $d$ of which are secretly using a new performance enhancing drug.  Furthermore, imagine that the only test for the drug is an expensive and time consuming blood test.  The trivial solution would be to collect blood samples from all $N$ athletes and then test each blood sample individually for the presence of the drug.  However, this is unnecessarily expensive when the test is accurate and the number of drug users is small.  A cheaper solution involves pooling portions of each player's blood into a small number of well-chosen testing pools.  Each of these testing pools can then be tested once, and the results used to identify the offenders.

A pooling-based testing procedure as described above can be modeled mathematically as a boolean matrix $\mathcal{M} \in \{ 0, 1 \}^{m \times N}$.  Each row of $\mathcal{M}$ corresponds to a subset of the $N$ athletes' whose blood will be pooled, mixed, and then tested once for the presence of the drug.  Hence, the goal of our nonadaptive group testing can be formulated at follows:  Design a matrix, $\mathcal{M} \in \{ 0, 1 \}^{m \times N}$, with as few rows as possible so that any boolean vector, $\vec{x} \in \{ 0, 1\}^N$, containing at most $d$ nonzero entries can be recovered exactly from the result of the pooled tests, $\mathcal{M} \vec{x} \in \{ 0,1 \}^m$.  Here all arithmetic is boolean, with the boolean $OR$ operator replacing summation and the boolean $AND$ operator replacing multiplication.  One well studied solution to this nonadaptive group testing problem is to let $\mathcal{M}$ be a $d$-disjunct matrix.

\begin{Definition}
An $m \times N$ binary matrix, $\mathcal{M} \in \{0,1\}^{m \times N}$, is called $d$-disjunct if for any subset of $d+1$ columns of $\mathcal{M}$, $C = \{ c_1,~c_2,~\dots,~c_{d+1} \} \subset [1,N] \cap \mathbbm{N}$, there exists a subset of $d+1$ rows of $\mathcal{M}$, $R = \{ j_1,~j_2,~\dots,~j_{d+1} \} \subset [1,m] \cap \mathbbm{N}$, such that the submatrix 
$$\left( \begin{array}{llll} 
\mathcal{M}_{j_1,c_1} & \mathcal{M}_{j_1,c_2} & \dots & \mathcal{M}_{j_1,c_{d+1}} \\ 
\mathcal{M}_{j_2,c_1} & \mathcal{M}_{j_2,c_2} & \dots & \mathcal{M}_{j_2,c_{d+1}} \\
& & \vdots & \\
\mathcal{M}_{j_{d+1},c_1} & \mathcal{M}_{j_{d+1},c_2} & \dots & \mathcal{M}_{j_{d+1},c_{d+1}} \\ 
\end{array} \right).$$
is the $(d+1) \times (d+1)$ identity matrix.\footnote{This is not the standard statement of the definition.  Traditionally, a boolean matrix $\mathcal{M}$ is said to be $d$-disjunct if the boolean $OR$ of any $d$ of its columns does not contain any other column \cite{FGroupTest,MyGroupTest}.  However, these two definitions are essentially equivalent.  The $d$-disjunct condition is also equivalent to the $\left(d + 1 \right)$-strongly selective condition utilized by compressed sensing algorithms based on group testing matrices \cite{MuthuCS}.}
\label{def:dDisjunct}
\end{Definition} 

Nonadaptive group testing is closely related to the recovery of ``exactly sparse'' vectors $\vec{x} \in \mathbbm{R}^N$ containing exactly $d$ nonzero entries.  In fact, it is not difficult to modify standard group testing techniques to solve such problems.  However, it is not generally possible to modify these approaches in order to obtain methods capable of achieving the type of approximation guarantees we are interested in here (i.e., see Equation~(\ref{eqn:Aerror2})).  However, fast $o(N)$-time approximation algorithms based on $d$-disjunct matrices with weaker approximation guarantees have been developed \cite{MuthuCS}.  Hence, if we can relate $\left( K,\alpha \right)$-coherent matrices to $d$-disjunct matrices, we will informally settle the design requirement regarding the existence of fast approximation algorithms (see the third design requirement in Section~\ref{sec:GenProb}).

\begin{lem}
An $m \times N$ $\left( K,\alpha \right)$-coherent matrix, $\mathcal{M}$, will also be $\lfloor (K - 1) / \alpha \rfloor$-disjunct.
\label{lem:CohDisjunct}
\end{lem}

\noindent \textit{Proof:}  Choose any subset of $\lfloor (K-1) / \alpha \rfloor+1$ columns from $\mathcal{M}$, $C = \{ c_1,~c_2,~\dots,~c_{\lfloor (K-1) / \alpha \rfloor+1} \} \subset [1,N] \cap \mathbbm{N}$.  Consider the column $\mathcal{M}_{\cdot,c_{1}} \in \{ 0,1 \}^m$.  Because $\mathcal{M}$ is a binary $\left( K,\alpha \right)$-coherent matrix, we know that there can be at most $\alpha$ rows, $j$, for which $\mathcal{M}_{j,c_1} ~=~ \mathcal{M}_{j,c_2}~=~1$.  Hence, there are at most $\alpha \lfloor (K - 1) / \alpha \rfloor ~\leq~ K - 1$ total rows in which $\mathcal{M}_{\cdot,c_1}$ will share a $1$ with any of the other columns listed in $C$.  Since $\mathcal{M}_{\cdot,c_{1}}$ contains at least $K$ ones, there exists a row, $j_1 \in [1,m] \cap \mathbbm{N}$, containing a $1$ in column $c_1$ and zeroes in all of $C - \{ c_1 \}$.  Repeating this argument with $c_2, \dots, c_{\lfloor (K-1) / \alpha \rfloor+1}$ replacing $c_1$ above proves the lemma.~~$\Box$ \\

Any $m \times N$ $d$-disjunct matrix must have $m = \Omega \left( \min \{ d^2 \log_d N, N \} \right)$ \cite{GTLowerBound}.  Furthermore, near-optimal explicit $d$-disjunct measurement matrix constructions of size $O(d^2 \log N) \times N$ exist \cite{ExplicitGT}.  Of more interest here, however, is that the lower bound for $d$-disjunct matrices together with Lemma~\ref{lem:CohDisjunct} provides a lower bound for $\left( K,\alpha \right)$-coherent matrices.  More specifically, we can see that any $m \times N$ $\left( K,\alpha \right)$-coherent matrix must have $m = \Omega \left( \min \left\{ (K^2 / \alpha^2) \log_{K / \alpha} N, N \right\} \right)$.

In the next section we will demonstrate that ideas from previous fast compressed sensing approximation methods based on $d$-disjunct matrices \cite{MuthuCS} can be utilized in combination with the properties of $\left( K,\alpha \right)$-coherent matrices to obtain the type of stronger approximation guarantees we consider in this paper.  In the process we will simultaneously decrease the previously obtained runtime complexities of these algorithms for general signals.  As a result, we will obtain entirely deterministic sublinear-time (in $N$) approximation algorithms which match the runtime and approximation guarantees previously only achieved with uniformly high probability by sublinear-time methods based on random measurement matrices (e.g., \cite{HHS}).

\subsection{Properties of Binary Incoherent Matrices}
\label{sec:relations}

The following theorem summarizes several important properties of $(K,\alpha)$-coherent matrices with respect to general sparse approximation problems.  Most importantly, the first statement guarantees the existence of a simple sublinear-time recovery algorithm, $\Delta_{\mathcal{M}}$, which is guaranteed to satisfy an approximation guarantee along the lines of Equation~\ref{eqn:Aerror2} for all $(K,\alpha)$-coherent matrices, $\mathcal{M}$, and vectors $\vec{x} \in \mathbbm{C}^N$.

\begin{Theorem}
Let $\mathcal{M}$ be an $m \times N$ $(K,\alpha)$-coherent matrix.  Then, all of the following statements will hold:
\begin{enumerate}

\item Let $\epsilon \in (0,1]$, $k \in \left[1,K \cdot \frac{\epsilon}{4 \alpha} \right) \cap \mathbbm{N}$.  There exists an approximation algorithm based on a modified form of $\mathcal{M}$, $\Delta_{\mathcal{M}}: \mathbbm{C}^{m \lceil \log_2 N \rceil + m} \rightarrow \mathbbm{C}^N$, that is guaranteed to output a vector $\vec{z}_{S} \in \mathbbm{C}^N$ satisfying
$$ \left\| \vec{x} - \vec{z}_{S} \right\|_2 ~\leq~ \left\| \vec{x} - \vec{x}^{\rm~ opt}_k \right\|_2 + \frac{22 \epsilon \left\| \vec{x} - \vec{x}^{\rm opt}_{(k/\epsilon)} \right\|_1}{\sqrt{k}}$$
for all $\vec{x} \in \mathbbm{C}^{N}$.  
Most importantly, $\Delta_{\mathcal{M}}$ can be evaluated in $O\left( m \log N \right)$-time.  See Appendix~\ref{appsec:SublinearProof} for details.

\item Define the $m \times N$ matrix $\mathcal{W}$ by normalizing the columns of $\mathcal{M}$ so that $\mathcal{W}_{i,j} = \mathcal{M}_{i,j} / \sqrt{\| \mathcal{M}_{\cdot,j} \|_1}$.  Then, the matrix $\mathcal{W}$ will be $\frac{\alpha}{K}$-coherent.

\item Furthermore, the $m \times N$ matrix $\mathcal{W}$ defined above will have the RIP$_{2}$($N$,$k$,$(k-1) \alpha / K$).\footnote{It is worth noting that modified $(K,\alpha)$-coherent matrices can also be used as Johnson-Lindenstrauss embeddings.  See \cite{RIPtoJL} together with \cite{JLtoRIP} to learn more about the near equivalence of Johnson-Lindenstrauss embeddings and RIP$_2$ matrices.}

\item Define the $m \times N$ matrix $\mathcal{W}$ by $\mathcal{W}_{i,j} = \mathcal{M}_{i,j} / \left(\| \mathcal{M}_{\cdot,j} \|_1 \right)^{\frac{1}{p}}$.  Then, the matrix $\mathcal{W}$ will have the\\ RIP$_{p}$($N$,$k$,$C(k-1)\alpha/K$) for all $1 \leq p \leq 1 + \frac{1}{\log N}$, where $C$ is an absolute constant larger than $1/2$.

\item $\mathcal{M}$ is $\lfloor (K - 1) / \alpha \rfloor$-disjunct.

\item $\mathcal{M}$ has at least $m = \Omega \left( \min \left\{ (K^2 / \alpha^2) \log_{K / \alpha} N, N \right\} \right)$ rows.

\end{enumerate}
\label{thm:MainResult}
\end{Theorem}

\noindent \textit{Proof:}  The proof of each part is as follows.
\begin{enumerate}
\item See Appendix~\ref{appsec:SublinearProof}.

\item The proof follows easily from the definitions.

\item The proof follows from part 2 together with Theorem~\ref{thm:RIPcoherent}.  However, for the sake of completeness we will recount the proof in more detail here.  

Let $X = \{ x_{1}, \dots, x_{k} \} \subset [0,N)$.  Given any such $X$ with $|X| = k$, we define $\mathcal{W}_{X}$ to be the $m \times k$ matrix consisting of the $k$ columns of $\mathcal{W}$ indexed by $X$.  We will consider the $k \times k$ Grammian (and therefore symmetric and non-negative definite) matrix
$$\mathcal{W}^T_{X} \mathcal{W}_{X} = \mathcal{I} + \mathcal{D}_{X}.$$  Our strategy will be to bound both $\| \mathcal{D}_{X} \|_{1}$ and $\| \mathcal{D}_{X} \|_{\infty}$ in the hope of applying Gerschgorin's theorem.

Each off diagonal entry $\left( \mathcal{D}_{X} \right)_{i,j}, i \neq j,$ is the inner product of $\mathcal{W}$'s $x_{i}$ and $x_{j}$ columns.  Thus, we have
$$\left( \mathcal{D}_{S} \right)_{i,j} = \frac{\mathcal{M}_{\cdot,i} \cdot \mathcal{M}_{\cdot,j}}{\sqrt{\| \mathcal{M}_{\cdot,i} \|_1 \| \mathcal{M}_{\cdot,j} \|_1}} \leq \frac{\alpha}{K}$$
since $\mathcal{M}$ is $(K,\alpha)$-coherent.  The end result is that both $\| \mathcal{D}_{X} \|_{1}$ and $\| \mathcal{D}_{X} \|_{\infty}$ are at most $\frac{(k-1) \cdot \alpha}{K}$.
Applying Gerschgorin's disk theorem we immediately see that the largest and smallest possible singular values of $\mathcal{W}_{X}$ are $\sqrt{1 + \frac{(k-1) \cdot \alpha}{K}}$ and $\sqrt{1 - \frac{(k-1) \cdot \alpha}{K}}$, respectively.  The result follows.

\item Note that we can consider $\mathcal{M}$ to be the adjacency matrix of a bipartite graph, $G = (A,B,E)$, with $|A| = N$ and $|B| = m$.  Each element of $A$ will have degree at least $K$.  Furthermore, for any $X \subset A$ with $|X| \leq k$ we can see that the set of neighbors of $X$ will have
$$|N(X)| \geq \sum^{|X|-1}_{j=0} ( K - j \alpha) \geq |X|\cdot K \cdot \left( 1 - \frac{\alpha(|X|-1)}{2K} \right).$$
Hence, $\mathcal{M}$ is the adjacency matrix of a $(k,K,(k-1)\alpha/2K)$-unbalanced expander graph.  The result now follows from the proof of Theorem 1 in \cite{RIP1expand}.
\end{enumerate}

\noindent Finally, the proof of parts 5 and 6 follow from Lemma~\ref{lem:CohDisjunct} and the subsequent discussion in Section~\ref{sec:GTcohMat}, respectively.~~$\Box$ \\

Recall that explicit constructions of $(K,\alpha)$-coherent matrices exist \cite{DetCS1,DetFRIP}.  It is worth noting that RIP$_2$ matrix constructions based on these $(K,\alpha)$-coherent matrices are optimal in the sense that any RIP$_2$ matrix with binary entries must have a similar number of rows \cite{BinaryBadRIP}.  More interestingly, Theorem~\ref{thm:MainResult} formally demonstrates that $(K,\alpha)$-coherent matrices satisfy all the Fourier design requirements in Section~\ref{sec:GenProb} other than the first one regarding small Fourier sampling requirements.   In the sections below we will consider an optimized number theoretic construction for $(K,\alpha)$-coherent matrices along the lines of the construction implicitly utilized in \cite{DetFRIP,ImpFourier}.  As we shall demonstrate, these constructions have small Fourier sampling requirements.  Hence, they will satisfy all four desired Fourier design requirements.

\section{A $(K,\alpha)$-Coherent Matrix Construction}
\label{sec:makeCoherent}

Let $\mathcal{F}_N$ denote the $N \times N$ unitary discrete Fourier transform matrix,
$$\left( \mathcal{F}_N \right)_{i,j} = \frac{\mathbbm{e}^{\frac{-2 \pi \mathbbm{i} \cdot i \cdot j}{N}}}{\sqrt{N}}.$$
Recall that we want an $m \times N$ matrix $\mathcal{M}$ with the property that $\mathcal{M} \mathcal{F}_N$ contains nonzero values in as few columns as possible.  In addition, we want $\mathcal{M}$ to be a binary $(K,\alpha)$-coherent matrix so that we can utilize the sublinear-time approximation technique provided by Theorem~\ref{thm:MainResult}.  It appears to be difficult to achieve both of these goals simultaneously as stated.  Hence, we will instead optimize a construction recently utilized in \cite{ImpFourier} which solves a trivial variant of this problem.

Let $\tilde{N}, N \in \mathbbm{N}$ with $\tilde{N} > N$.  We will say that an $m \times \tilde{N}$ matrix, $\widetilde{\mathcal{M}}$, is $(K,\alpha)_N$-coherent if the $m \times N$ submatrix of $\widetilde{\mathcal{M}}$ formed by its first $N$ columns is $(K,\alpha)$-coherent.  In what follows we will consider ourselves to be working with $(K,\alpha)_N$-coherent matrices whose first $N$ rows match a given $m \times N$ $(K,\alpha)$-coherent matrix, $\mathcal{M}$, of interest.  Note that this slight generalization will not meaningfully change anything previously discussed.  For example, we may apply Theorem~\ref{thm:MainResult} to the submatrix formed by the first $N$ columns of any given $(K,\alpha)_N$-coherent matrix, $\widetilde{\mathcal{M}}$, thereby effectively applying Theorem~\ref{thm:MainResult} to $\widetilde{\mathcal{M}}$ in the context of approximating vectors belonging to a fixed $N$-dimensional subspace of $\mathbbm{C}^{\tilde{N}}$.  The last $\tilde{N} - N$ columns of any $(K,\alpha)_N$-coherent matrix $\widetilde{\mathcal{M}}$ will be entirely ignored throughout this paper with one exception:  We will hereafter consider it sufficient to guarantee that $\widetilde{\mathcal{M}} \mathcal{F}_{\tilde{N}}$ (as opposed to $\mathcal{M} \mathcal{F}_N$) contains nonzero values in as few columns as possible.  This modification will not alter the sparse Fourier approximation guarantees (i.e., see Equation~(\ref{eqn:Aerror})) obtainable via Theorem~\ref{thm:MainResult} in any way when the functions being approximated are $N$-bandlimited.  However, allowing $\tilde{N}$ to be greater than $N$ will help us obtain small Fourier sampling requirements.

Let $\widetilde{\mathcal{M}}$ be an $m \times \tilde{N}$ $(K,\alpha)_N$-coherent matrix.  It is useful to note that the column sparsity we desire in $\widetilde{\mathcal{M}} \mathcal{F}_{\tilde{N}}$ is closely related to the discrete uncertainty principles previously considered in \cite{DonStark}.

\begin{Theorem} \textrm{\textbf{(See \cite{DonStark}).}}
Suppose $\vec{y} \in \mathbbm{C}^{\tilde{N}}$ contains $\tilde{N}_t$ nonzero entries, while $\hat{y} = \vec{y}^{~\rm T} \mathcal{F}_{\tilde{N}}$ contains $\tilde{N}_{\omega}$ nonzero entries.  Then, $\tilde{N}_t \tilde{N}_{\omega} \geq \tilde{N}$.  Furthermore, $\tilde{N}_t \tilde{N}_{\omega} = \tilde{N}$ holds if and only if $\vec{y}$ is a scalar multiple of a cyclic permutation of the picket fence sequence in $\mathbbm{C}^{\tilde{N}}$ containing $v$ equally-spaced nonzero elements 
$$\left( III^v \right)_u ~=~\left\{ \begin{array}{ll} 1 & \textrm{if } u \equiv 0 \textrm{ mod } \frac{\tilde{N}}{v} \\ 0 & {\rm otherwise} \end{array} \right., $$
where $v \in \mathbbm{N}$ divides $\tilde{N}$.
\label{thm:FourierUP}
\end{Theorem}

We will build $m \times \tilde{N}$ $(K,\alpha)_N$-coherent matrices, $\mathcal{M}$, below whose rows are each a permuted binary picket fence sequence.  In this case Theorem~\ref{thm:FourierUP} can be used to bound the number of columns of $\widetilde{\mathcal{M}} \mathcal{F}_{\tilde{N}}$ which contain nonzero entries.  This, in turn, will bound the number of function samples required in order to approximate a given periodic bandlimited function.

We create an $m \times \tilde{N}$ $(K,\alpha)_N$-coherent matrix $\mathcal{M}$ as follows:  Choose $K$ pairwise relatively primes integers
$$s_1 < \cdots < s_K$$
and let $\tilde{N} = \prod^{K}_{j = 1} s_j ~>~ N$.
Next, we produce a picket fence row, $r_{j,h}$, for each $j \in [1,K] \cap \mathbbm{N}$ and $h \in [0,s_j) \cap \mathbbm{Z}$.  Thus, the $n^{\rm th}$ entry of each row $r_{j,h}$ is given by
\begin{equation}
(r_{j,h})_{n} = \delta \left( (n-h) \textrm{ mod } s_j \right) = \left\{ \begin{array}{ll} 1 & \textrm{if } n \equiv h \textrm{ mod } s_j \\ 0 & {\rm otherwise} \end{array} \right.,
\label{eqn:Def_r}
\end{equation}
where $n \in [0,\tilde{N}) \cap \mathbbm{Z}$. 
We then form $\mathcal{M}$ by setting
\begin{equation}
\mathcal{M} = \left( \begin{array}{l} r_{1,0} \\ r_{1,1}\\ \vdots \\ r_{1,s_{1}-1} \\ r_{2,0} \\ \vdots \\ r_{2,s_{2}-1} \\ \vdots \\ r_{K,s_{K}-1} \\
\end{array} \right).
\label{eqn:Def_M}
\end{equation}
For an example measurement matrix see Figure~\ref{fig:RIPexample}.

\begin{figure}
$$\textbf{---------------------------------------------------------------------------}$$
$$\begin{array}{llllllllll}
~\mathbf{n} \in \mathbf{[0,\tilde{N})}&\hspace{16pt}& \mathbf{0} & \mathbf{1} & \mathbf{2} & \mathbf{3} & \mathbf{4} & \mathbf{5} & \mathbf{6} & \dots \\
\end{array}$$
$$\begin{array}{l} \mathbf{n} \equiv \mathbf{0}~(\mathbf{mod}~\mathbf{2}) \\ \mathbf{n} \equiv \mathbf{1}~(\mathbf{mod}~\mathbf{2}) \\ \mathbf{n} \equiv \mathbf{0}~(\mathbf{mod}~\mathbf{3}) \\ \mathbf{n} \equiv \mathbf{1}~(\mathbf{mod}~3) \\ \mathbf{n} \equiv \mathbf{2}~(\mathbf{mod}~\mathbf{3}) \\ \quad \vdots \\ \mathbf{n} \equiv \mathbf{1}~(\mathbf{mod}~\mathbf{5}) \\ \quad \vdots \end{array}
\begin{array}{l}  \\  \\  \\  \\  \\  \\  \\  \end{array}
\left( \begin{array}{llllllll}
1 & 0 & 1 & 0 & 1 & 0 & 1 & \dots \\
0 & 1 & 0 & 1 & 0 & 1 & 0 & \dots \\
1 & 0 & 0 & 1 & 0 & 0 & 1 & \dots \\
0 & 1 & 0 & 0 & 1 & 0 & 0 & \dots \\
0 & 0 & 1 & 0 & 0 & 1 & 0 & \dots \\
&&& \vdots &&&& \\
0 & 1 & 0 & 0 & 0 & 0 & 1 & \dots \\
&&& \vdots &&&& \\
\end{array} \right)$$
\caption{Measurement Matrix, $\mathcal{M}$, Using $s_{1} = 2$, $s_{2} = 3$, $s_{3} = 5$, \dots}
\label{fig:RIPexample}
$$\textbf{---------------------------------------------------------------------------}$$
\end{figure}

\begin{lem}
An $m \times \tilde{N}$ matrix $\mathcal{M}$ as constructed in Equation~(\ref{eqn:Def_M}) will be $(K, \lfloor \log_{s_1} N \rfloor)_N$-coherent with $m = \sum^K_{j = 1} s_j$.
\label{lem:NumCon1}
\end{lem}

\noindent \textit{Proof:}  Choose any two distinct integers, $l \neq n$, from $[0,N)$.  Let $\mathcal{M}_{\cdot,l}$ and $\mathcal{M}_{\cdot,n}$ denote the $l^{\rm th}$ and $n^{\rm th}$ columns of $\mathcal{M}$, respectively.  The inner product of these columns is
$$\mathcal{M}_{\cdot,l} \cdot \mathcal{M}_{\cdot,n} = \sum^K_{j = 1} \delta \left( (n-l)~{\rm mod}~s_j \right).$$
The sum above is at most the maximum $\alpha$ for which $\prod^\alpha_{j=1} s_j \leq N$ by the Chinese Remainder Theorem.  Furthermore, this value is itself bounded above by $\lfloor \log_{s_1} N \rfloor$.  The equation for $m$ immediately follows from the construction of $\mathcal{M}$ above.~~$\Box$ \\

The following Lemma is a consequence of Theorem~\ref{thm:FourierUP}.  

\begin{lem}
Let $\mathcal{M}$ be an $m \times \tilde{N}$ matrix as constructed in Equation~(\ref{eqn:Def_M}).  Then, $\mathcal{M} \mathcal{F}_{\tilde{N}}$ will contain nonzero entries in exactly $m - K + 1 = \left( \sum^K_{j = 1} s_j \right) - K + 1$ columns.
\label{lem:SmallFourier}
\end{lem}

\noindent \textit{Proof:}  Fix $j \in [1,K] \cap \mathbbm{N}$. Each picket fence row, $r_{j,h} \in \{ 0, 1 \}^{\tilde{N}}$, contains $\tilde{N} / s_j$ ones.  Thus, $r_{j,h}^{\rm T} \mathcal{F}_{\tilde{N}}$ contains $s_j$ nonzero entries for all $h \in [0,s_j) \cap \mathbbm{Z}$.  Furthermore, $r_{j,h}^{\rm T} \mathcal{F}_{\tilde{N}}$ contains nonzero values \textit{in the same entries for all $h \in [0,s_j) \cap \mathbbm{Z}$} since all $r_{j,h}$ rows (with $j$ fixed) are cyclic permutations of one another.  Finally, let $l,j \in [1,K] \cap \mathbbm{N}$ with $j \neq l$ and suppose that $r_{j,h}^{\rm T} \mathcal{F}_{\tilde{N}}$ and $r_{l,g}^{\rm T} \mathcal{F}_{\tilde{N}}$ both have nonzero values in the same entry.  This can only happen if
$$h \frac{\tilde{N}}{s_j} = g \frac{\tilde{N}}{s_l}$$
for a pair of integers $0 \leq h < s_j$ and $0 \leq g < s_l$.
However, since $s_j$ and $s_l$ are relatively prime, Euclid's lemma implies that this can only happen when $h = g = 0$.  The result follows.~~$\Box$ \\

We can now see that matrix construction presented in this section satisfies all four of our Fourier design requirements.  In the next sections we will consider methods for optimizing the relatively prime integer values, $s_1, \dots, s_K$, used to construct our $(K,\alpha)_N$-coherent matrices.  In what follows we will drop the slight distinction between $(K,\alpha)_N$-coherent and $(K,\alpha)$-coherent matrices for ease of discussion. 

\section{Optimizing the $(K,\alpha)$-Coherent Matrix Construction}
\label{sec:OptBuild}

Note that $K$ appears as part of a ratio involving $\alpha$ in each statement of Theorem~\ref{thm:MainResult}.  Hence, we will focus on constructing $(K,\alpha)$-coherent matrices in which $K$ is a constant multiple of $\alpha$ in this section.  For a given value of $D \in (1, \infty)$ we can optimize the Section~\ref{sec:makeCoherent} methods for constructing a $\left(D\alpha,\alpha \right)$-coherent matrix with a small number of rows by reformulating the matrix design problem as an optimization problem (see Figure~\ref{fig:Problem}).  In this section we will develop concrete bounds for the number of rows, $m$ as a function of $D$, $N$, and $\alpha$, that will appear in any $m \times N$ $\left(K = D\alpha,\alpha \right)$-coherent matrix constructed as per Section~\ref{sec:makeCoherent}.  These bounds will ultimately allow us to cast the matrix optimization problem in Figure~\ref{fig:Problem} as a linear integer program in Section~\ref{sec:LinIntProg}.

\begin{figure}
$$\textbf{---------------------------------------------------------------------------}$$
\centering
Minimize
\begin{equation}
m = \sum^{K_{\alpha} = \left \lceil D \alpha \right \rceil}_{j = 1} s_j
\label{eqn:Optsamps}
\end{equation}
subject to the following constraints:\\
\begin{enumerate}
\renewcommand{\theenumi}{\Roman{enumi}}
\centering
%\item \hspace{.9 in} $\alpha \in \left\{ 1, 2, \dots, \lfloor \log_2 N \rfloor \right\}$.\\
%\item \hspace{1.45 in} Set $\frac{(k-1)\alpha}{\epsilon}$.\\
\item \hspace{1.45 in} $s_1 < \cdots < s_{K}$.\\
\item \hspace{.85 in} $\prod^{\alpha}_{j=1} s_j < N \leq \prod^{\alpha+1}_{j=1} s_j$.\\
\item $s_1, \cdots, s_{K}$ are pairwise relatively prime.\\
\end{enumerate}
$$\textbf{---------------------------------------------------------------------------}$$
\caption{Matrix Design Optimization Problem for Given $N$, $D$, and $\alpha$ values}
\label{fig:Problem}
\end{figure}

The following trivial fact will be useful below.

\begin{lem}
Let $x_1, x_2, \dots, x_n \in [2,\infty)$ be such that $x_n \geq x_{n-1} \geq \dots \geq x_1 \geq 2$.  Then, $\sum^n_{j=1} x_j \leq \prod^n_{j=1} x_j$.
\label{lem:SimpBound}
\end{lem}

\noindent \textit{Proof:}  This follows immediately from the fact that
$$1 + \frac{\sum^{n-1}_{j=1} x_j}{x_n} ~\leq~ n ~\leq~ 2^{n-1} ~\leq~ \prod^{n-1}_{j=1} x_j. \hspace{0.5 in} \Box$$

Define $p_{0} = 1$ and let $p_{l}$ be the $l^{\rm th}$ prime natural number.  Thus, we have
\begin{equation}
p_{0} = 1, p_{1} = 2, p_{2} = 3, p_{3} = 5, p_{4} = 7, \dots
\label{eqn:Primes}
\end{equation}
Suppose that $S = \{ s_1, \dots, s_{K} \}$ is a solution to the optimization problem presented in Figure~\ref{fig:Problem} for given values of $\alpha$, $D$, and $N$.  Let $p_{q_S}$ be the largest prime factor appearing in any element of $S$.  Finally, let
$$q = \max \left\{ q_{S} ~~\big|~~ S \textrm{ solves the optimization problem in Figure~\ref{fig:Problem}} \right\}.$$
The following lemma bounds $q$ as a function of $N$, $D$, and $\alpha$.

\begin{lem}
Suppose $s_1, s_2, \cdots, s_{K}$ satisfy all three constraints in Figure~\ref{fig:Problem}.  Set $\tilde{m} = \sum^{K}_{j=1} s_{j}$.  Next, let $p_t$ be the smallest prime number greater than $2$ for which
$$p_t \cdot \left( K - \alpha - 1 \right) + \left( K - \alpha - 1 \right)\left( K - \alpha - 2 \right) + \left( \alpha + 1 \right) N^{\frac{1}{\alpha+1}} > \tilde{m}$$
holds.  Then, $q < t + K - \alpha - 1$.
\label{lem:PrimeUpperBound}
\end{lem}

\noindent \textit{Proof:}\\

Let $s'_1, s'_2, \cdots, s'_{K}$ be a solution to the optimization problem presented in Figure~\ref{fig:Problem}. Set $m = \sum^{K}_{j=1} s'_{j}$.  Note that there must exist at least one prime, $p_{\tilde{l}} \in [p_t,p_{t + K - \alpha - 1})$, which is not a prime factor of any $s'_{j}$ value.  If no such prime exists, then Lemma~\ref{lem:SimpBound} applied to the prime factors of each $s'_j$ containing one of the primes in $[p_t,p_{t + K - \alpha - 1})$ tells us that the sum of $s'_1, s'_2, \cdots, s'_{K}$ must be
$$m ~\geq~ \sum^{K - \alpha - 2}_{j=0} p_{t+j} ~+~ \sum^{\alpha+1}_{j=1} s'_j.$$
The second constraint in Figure~\ref{fig:Problem} together with the arithmetic-geometric mean inequality tells that we must always have
\begin{equation}
\left( \alpha + 1 \right) N^{\frac{1}{\alpha+1}} ~\leq~ \left( \alpha + 1 \right) \cdot \left( \prod^{\alpha+1}_{j=1} s'_j \right)^{\frac{1}{\alpha+1}} ~\leq~ \sum^{\alpha+1}_{j=1} s'_j.
\label{eqn:HeadBound}
\end{equation}
Furthermore, it is not difficult to see that
\begin{equation}
\label{eqn:TailBound0}
\sum^{K - \alpha - 2}_{j=0} p_{t+j} ~\geq~ \sum^{K - \alpha - 2}_{j=0} \left( p_{t} + 2j \right) ~\geq~  p_{t} \cdot \left(K - \alpha - 1 \right) ~+~ \left( K - \alpha - 1 \right)\left(K - \alpha - 2\right)
\end{equation}
since $p_t > 2$.  Thus, if every prime in $[p_t,p_{t + K - \alpha - 1})$ appears as a prime factor in some $s'_{j}$, then
$$m ~\geq~ p_t \cdot \left( K - \alpha - 1 \right) + \left( K - \alpha - 1 \right)\left( K - \alpha - 2 \right) + \left( \alpha + 1 \right) N^{\frac{1}{\alpha+1}} ~>~ \tilde{m},$$
violating our assumption concerning the optimality of $s'_1, s'_2, \cdots, s'_{K}$.  This proves our claim regarding the existence of at least one prime, $p_{\tilde{l}} \in [p_t,p_{t + K - \alpha - 1})$, which is not a prime factor of any $s'_{j}$ value.

Now suppose that some $s'_{j'}$ contains a prime factor, $p_{l'}$, with $l' \geq t + K - \alpha - 1$.  Substitute the largest currently unused prime, $p_{\tilde{l}} \in [p_t,p_{t + K - \alpha - 1})$, for $p_{l'}$ in the prime factorization of $s'_{j'}$ to obtain a smaller value, $s'_{\tilde{j}}$.  If we can show that $s'_1, s'_2, \cdots, s'_{K}$ with $s'_{\tilde{j}}$ substituted for $s'_{j'}$ still satisfies all three Figure~\ref{fig:Problem} constraints after reordering, we will again have a contradiction to the assumed minimality of our original solution.  In fact, it is not difficult to see that all constraints other than II above will trivially be satisfied by construction.  Furthermore, if $s'_{\tilde{j}} > s'_{\alpha+1}$, then Constraint II will also remain satisfied and we will violate our assumption that the $s'_j$ values originally had a minimal sum.

Finally, the second case where $p_t \leq p_{\tilde{l}} \leq s'_{\tilde{j}} < s'_{\alpha+1}$ could only occur if originally
\begin{equation}
\sum^{K}_{j=\alpha+2} s'_{j} \geq \sum^{K - \alpha - 1}_{j=1} \left( p_t + 2j \right) = p_t \cdot \left( K - \alpha - 1 \right) + \left( K - \alpha\right)\left( K - \alpha - 1 \right).
\label{eqn:TailBound}
\end{equation}
When combined with Equation~(\ref{eqn:HeadBound}) above, Equation~(\ref{eqn:TailBound}) reveals that if $s'_{\tilde{j}} < s'_{\alpha+1}$ then we must have originally had
$$\sum^{K}_{j=1} s'_j \geq \left( \alpha + 1 \right) N^{\frac{1}{\alpha+1}} + p_t \cdot \left( K - \alpha - 1 \right) + \left( K - \alpha\right)\left( K - \alpha - 1 \right) > \tilde{m}.$$
However, in this case the assumed minimality of $s'_1, s'_2, \dots, s'_{K}$ would again have been violated.~~$\Box$ \\

We will now establish a slightly more refined result than that of Lemma \ref{lem:PrimeUpperBound}.

\begin{lem}
\label{lem:PrimeUpperBound2}
Suppose $s_1, s_2, \cdots, s_{K}$ satisfy all three constraints in Figure~\ref{fig:Problem}.  Set $\tilde{m} = \sum^{K}_{j=1} s_{j}$.  Let $L=\prod_{i=1}^v p_i$ for any desired $v \in \mathbb{N}$, and let $\phi(L) = \prod_{i=1}^v (p_i-1).$
Next, let $p_t$ be the smallest prime number greater than 2 for which
\begin{align*}
p_t \cdot \left(K - \alpha - 1 \right)
+
\left( K - \alpha -1 \right)\left(K - \alpha -2 \right)
+
 \frac{\big(L - 2 \phi(L) -2v \big) (\phi(L)+v)
\big\lfloor \frac{K - \alpha - 2}{\phi(L)+v} \big\rfloor
\Big(\big\lfloor \frac{K - \alpha - 2}{\phi(L)+v} \big\rfloor  -1 \Big)
  }{2} \qquad
&
\\
 \qquad \qquad +
\big(L - 2 \phi(L)-2v \big) \left\lfloor \frac{K - \alpha - 2}{\phi(L)+v} \right\rfloor \left( K-\alpha -1 - (\phi(L)+v) \left\lfloor \frac{K - \alpha - 2}{\phi(L)+v} \right\rfloor \right) + \left( \alpha + 1 \right) N^{\frac{1}{\alpha+1}} & > \tilde{m}
\end{align*}
holds.  Then, $q < t + K - \alpha - 1$.
\end{lem}

\noindent \textit{Proof:}\\

We will prove this lemma by modifying our proof of Lemma \ref{lem:PrimeUpperBound}.  In particular, we will modify formulas (\ref{eqn:TailBound0}) and (\ref{eqn:TailBound}).
Note that amongst any $L$ consecutive numbers, there are at most $\phi(L)+v$ prime numbers.  Hence, we have that $p_{i+w (\phi(L)+v)} \ge p_{i} + wL$ for all $w \in \mathbbm{N}$.  Thus, we may replace formula (\ref{eqn:TailBound0}) with
\begin{align*}
\sum^{K - \alpha - 2}_{j=0} p_{t+j}
&\geq
\sum^{K - \alpha - 2}_{j=0} \left(p_{t+j- (\phi(L)+v) \left\lfloor \frac{j}{\phi(L)+v} \right\rfloor } + L \left\lfloor \frac{j}{\phi(L)+v} \right\rfloor \right)
\\ &
\geq
\sum^{K - \alpha - 2}_{j=0} \left(p_{t} + 2 \left(j- (\phi(L)+v) \left\lfloor \frac{j}{\phi(L)+v} \right \rfloor \right) + L \left\lfloor \frac{j}{\phi(L)+v} \right\rfloor \right)
\\ &
= \sum^{K - \alpha - 2}_{j=0} \left(p_{t} + 2j +\big(L - 2 \phi(L) -2v \big) \left\lfloor \frac{j}{\phi(L)+v} \right\rfloor  \right)
\\&
= p_t \cdot \left(K - \alpha - 1 \right) + \left( K - \alpha - 1 \right)\left(K - \alpha - 2\right) + \big(L - 2 \phi(L) -2v \big) \sum^{K - \alpha - 2}_{j=0} \left\lfloor \frac{j}{\phi(L)+v} \right\rfloor
\\&
= p_t \cdot \left(K - \alpha - 1 \right)
+
\left( K - \alpha -1 \right)\left(K - \alpha -2 \right)
+
 \frac{\big(L - 2 \phi(L) -2v \big) (\phi(L)+v)
\big\lfloor \frac{K - \alpha - 2}{\phi(L)+v} \big\rfloor
\Big(\big\lfloor \frac{K - \alpha - 2}{\phi(L)+v} \big\rfloor  -1 \Big)
  }{2}
\\
& \qquad \qquad +
\big(L - 2 \phi(L)-2v \big) \left\lfloor \frac{K - \alpha - 2}{\phi(L)+v} \right\rfloor \left( K-\alpha -1 - (\phi(L)+v) \left\lfloor \frac{K - \alpha - 2}{\phi(L)+v} \right\rfloor \right).
\end{align*}
Note that amongst any $L$ consecutive numbers, a maximal subset of pairwise relatively prime numbers has at most $\phi(L)+v$ numbers.
Hence, we may also replace formula (\ref{eqn:TailBound}) by a similar argument to above with
\begin{align*}
\sum^{K}_{j=\alpha+2} s'_{j}
& =
\sum^{K - \alpha - 2}_{j=0} s'_{\alpha+2+j}
\\ & \geq
\sum^{K - \alpha - 2}_{j=0} \left(s'_{\alpha+2+j- (\phi(L)+v) \left\lfloor \frac{j}{\phi(L)+v} \right\rfloor } + L \left\lfloor \frac{j}{\phi(L)+v} \right\rfloor \right)
\\ &
\geq
\sum^{K - \alpha - 2}_{j=0} \left(p_t + 2 \left(1+ j- (\phi(L)+v) \left\lfloor \frac{j}{\phi(L)+v} \right \rfloor \right) + L \left\lfloor \frac{j}{\phi(L)+v} \right\rfloor \right)
\\ &
= \sum^{K - \alpha - 2}_{j=0} \left(p_{t} + 2(1+j) +\big(L - 2 \phi(L) -2v \big) \left\lfloor \frac{j}{\phi(L)+v} \right \rfloor  \right)
\\&
= p_t \cdot \left(K - \alpha - 1 \right) + \left( K - \alpha \right)\left(K - \alpha - 1\right) + \big(L - 2 \phi(L) -2v \big) \sum^{K - \alpha - 2}_{j=0} \left\lfloor \frac{j}{\phi(L)+v} \right\rfloor
\\&
= p_t \cdot \left(K - \alpha - 1 \right)
+
\left( K - \alpha \right)\left(K - \alpha - 1\right)
+
 \frac{\big(L - 2 \phi(L) -2v \big) (\phi(L)+v)
\big\lfloor \frac{K - \alpha - 2}{\phi(L)+v} \big\rfloor
\Big(\big\lfloor \frac{K - \alpha - 2}{\phi(L)+v} \big\rfloor  -1 \Big)
  }{2}
\\
& \qquad \qquad +
\big(L - 2 \phi(L)-2v \big) \left\lfloor \frac{K - \alpha - 2}{\phi(L)+v} \right\rfloor \left( K-\alpha -1 - (\phi(L)+v) \left\lfloor \frac{K - \alpha - 2}{\phi(L)+v} \right\rfloor \right).
\end{align*}
By replacing (\ref{eqn:TailBound0}) and (\ref{eqn:TailBound}) with these bounds in the proof of Lemma \ref{lem:PrimeUpperBound}, we obtain the desired result.~~$\Box$ \\

The following corollary of Lemma~\ref{lem:PrimeUpperBound} provides a simple initial upper bound on the largest prime factor that may appear in any solution to the optimization problem presented in Figure~\ref{fig:Problem}.

\begin{cor}
Let $r$ be such that
$\prod^{\alpha}_{j=1} p_{r+j} < N \leq \prod^{\alpha+1}_{j=1} p_{r+j},$
and set $\tilde{m} = \sum^{K}_{j=1} p_{r+j}$.  Next, let $p_t$ be the smallest prime larger than $2$ for which
$$p_t \cdot \left( K - \alpha - 1 \right) + \left( K - \alpha - 1 \right)\left( K - \alpha - 2 \right) + \left( \alpha + 1 \right) N^{\frac{1}{\alpha+1}} > \tilde{m}$$
holds.  Then, $q < t + K - \alpha - 1$.
\label{cor:PrimeUpperBound}
\end{cor}

\noindent \textit{Proof:}\\

It is not difficult to see that
\begin{equation}
s_1 = p_{r+1}, s_2 = p_{r+2}, \cdots, s_{K} = p_{r+K}
\label{eqn:PrimeSoln}
\end{equation}
collectively satisfy all three constraints in Figure~\ref{fig:Problem}.  Applying Lemma~\ref{lem:PrimeUpperBound} yields the stated result.~~$\Box$ \\

Similarly, one can obtain the following corollary from Lemma \ref{lem:PrimeUpperBound2}.

\begin{cor}
\label{cor:PrimeUpperBound2}
Let $r$ be such that
$\prod^{\alpha}_{j=1} p_{r+j} < N \leq \prod^{\alpha+1}_{j=1} p_{r+j},$
and set $\tilde{m} = \sum^{K}_{j=1} p_{r+j}$.  Let $L=\prod_{i=1}^v p_i$ for any $v \in \mathbb{N}$, and let $\phi(L) = \prod_{i=1}^v (p_i-1).$  Next, let $p_t$ be the smallest prime larger than $2$ for which
\begin{align*}
p_t \cdot \left(K - \alpha - 1 \right)
+
\left( K - \alpha -1 \right)\left(K - \alpha -2 \right)
+
 \frac{\big(L - 2 \phi(L) -2v \big) (\phi(L)+v)
\big\lfloor \frac{K - \alpha - 2}{\phi(L)+v} \big\rfloor
\Big(\big\lfloor \frac{K - \alpha - 2}{\phi(L)+v} \big\rfloor  -1 \Big)
  }{2} \qquad
&
\\
 \qquad \qquad +
\big(L - 2 \phi(L)-2v \big) \left\lfloor \frac{K - \alpha - 2}{\phi(L)+v} \right\rfloor \left( K-\alpha -1 - (\phi(L)+v) \left\lfloor \frac{K - \alpha - 2}{\phi(L)+v} \right\rfloor \right) + \left( \alpha + 1 \right) N^{\frac{1}{\alpha+1}} & > \tilde{m}
\end{align*}
holds.  Then, $q < t + K - \alpha - 1$.
\end{cor}

The following lemma provides upper and lower bounds for the members of any valid solution to the optimization problem in Figure~\ref{fig:Problem} as functions of $N$, $D$, and $\alpha$.  This lemma is critical to the formulation of the optimization problem in Figure~\ref{fig:Problem} as a linear integer program in Section~\ref{sec:LinIntProg}.

\begin{lem}
The following bounds hold for any valid solution, $S = \{ s_1, s_2, \dots, s_K \}$, to the optimization problem in Figure~\ref{fig:Problem}:
\begin{enumerate}
\item $s_1 < \cdots < s_{K}$.
\item $s_1 \geq 2, s_2 \geq 3, \dots, s_K \geq p_K$.
\item $s_1 < N^{\frac{1}{\alpha}}$ and $s_{\alpha+1} > N^{\frac{1}{\alpha+1}}$.
\item Let $t \in \mathbbm{N}$ be defined as in Lemma~\ref{lem:PrimeUpperBound}, Lemma~\ref{lem:PrimeUpperBound2}, Corollary~\ref{cor:PrimeUpperBound}, or Corollary~\ref{cor:PrimeUpperBound2}.  Then, $s_K < p_{t + K - \alpha - 1}$.
\end{enumerate}
\label{lem:sjBounds}
\end{lem}

\noindent \textit{Proof:}\\

Assertion $(1)$ is a restatement of Constraint I in Figure~\ref{fig:Problem}.  The second assertion follows immediately from the fact that the ordered $s_j$ values must be pairwise relatively prime (i.e., Constraint III).  The third assertion follows easily from Constraint II.  Assertion $(4)$ follows from an argument analogous to the proof of Lemma~\ref{lem:PrimeUpperBound}.  That is, if $s_K \geq p_{t + K - \alpha - 1}$, then we may substitute $s_K$ with the largest prime in $[p_t,p_{t + K - \alpha - 1})$ which is not currently a prime factor of $s_1, \dots, s_K$ and thereby derive a contradiction.~~$\Box$ \\

The following lemmas provide concrete lower bound for $m$ in terms of $N$, $D$, and $\alpha$ (see Equation~(\ref{eqn:Optsamps}) in Figure~\ref{fig:Problem}).  These lemmas will ultimately allow us to judge the possible performance of any solution to our optimization problem based solely on the value of $\alpha$ whenever $N$ and $D$ are fixed.

\begin{lem}
Any solution to the optimization problem in Figure~\ref{fig:Problem} must have
$$m ~\geq~ K N^{\frac{1}{\alpha+1}} + \left( K - \alpha \right) \left( K - \alpha -1 \right).$$
\label{lem:mLowerBound}
\end{lem}

\noindent \textit{Proof:}\\

We know that $s_{\alpha + 2} > s_{\alpha+1} > N^{\frac{1}{\alpha+1}}$ from Lemma~\ref{lem:sjBounds}.  Hence, we can see that
$$\sum^{K}_{j= \alpha + 2} s_j ~\geq~ \sum^{K - \alpha - 1}_{j=1} \left( N^{\frac{1}{\alpha+1}} + 2j \right) \geq \left( K - \alpha - 1 \right) N^{\frac{1}{\alpha+1}} + \left( K - \alpha \right) \left( K - \alpha -1 \right).$$
Combining this lower bound with Equation~(\ref{eqn:HeadBound}) proves the lemma.~~$\Box$ \\

\begin{cor}  Let $L=\prod_{i=1}^v p_i$ for any desired $v \in \mathbb{N}$, and let $\phi(L) = \prod_{i=1}^v (p_i-1).$  Any solution to the optimization problem in Figure~\ref{fig:Problem} must have
\begin{align*}
m & ~>~  K N^{\frac{1}{\alpha+1}} 
+
\left( K - \alpha \right)\left(K - \alpha - 1\right)
+
 \frac{\big(L - 2 \phi(L) -2v \big) (\phi(L)+v)
\big\lfloor \frac{K - \alpha - 2}{\phi(L)+v} \big\rfloor
\Big(\big\lfloor \frac{K - \alpha - 2}{\phi(L)+v} \big\rfloor  -1 \Big)
  }{2}
\\
& \qquad \qquad +
\big(L - 2 \phi(L)-2v \big) \left\lfloor \frac{K - \alpha - 2}{\phi(L)+v} \right\rfloor \left( K-\alpha -1 - (\phi(L)+v) \left\lfloor \frac{K - \alpha - 2}{\phi(L)+v} \right\rfloor \right).
\end{align*}
\label{cor:mLowerBound2}
\end{cor}

\noindent \textit{Proof:}\\

Note that amongst any $L$ consecutive integers, a maximal subset of pairwise relatively prime numbers has at most $\phi(L)+v$ numbers.  As in the proof of Lemma \ref{lem:mLowerBound}, this corollary follows by combining Equation~(\ref{eqn:HeadBound}) with the fact that
\begin{align*}
\sum^{K}_{j= \alpha + 2} s_j 
&  =
\sum^{K - \alpha - 2}_{j=0} s_{\alpha+2+j}
\\ & \geq
\sum^{K - \alpha - 2}_{j=0} \left(s_{\alpha+2+j- (\phi(L)+v) \left\lfloor \frac{j}{\phi(L)+v} \right\rfloor } + L \left\lfloor \frac{j}{\phi(L)+v} \right\rfloor \right)
\\ &
\geq
\sum^{K - \alpha - 2}_{j=0} \left( s_{\alpha+1} + 2 \left(1+ j- (\phi(L)+v) \left\lfloor \frac{j}{\phi(L)+v} \right \rfloor \right) + L \left\lfloor \frac{j}{\phi(L)+v} \right\rfloor \right)
\\ &
> \sum^{K - \alpha - 2}_{j=0} \left( N^{\frac{1}{\alpha+1}}+ 2(1+j) +\big(L - 2 \phi(L) -2v \big) \left\lfloor \frac{j}{\phi(L)+v} \right \rfloor  \right)
\\&
=  N^{\frac{1}{\alpha+1}} \cdot \left(K - \alpha - 1 \right) + \left( K - \alpha \right)\left(K - \alpha - 1\right) + \big(L - 2 \phi(L) -2v \big) \sum^{K - \alpha - 2}_{j=0} \left\lfloor \frac{j}{\phi(L)+v} \right\rfloor
\\&
=  N^{\frac{1}{\alpha+1}} \cdot \left(K - \alpha - 1 \right)
+
\left( K - \alpha \right)\left(K - \alpha - 1\right)
+
 \frac{\big(L - 2 \phi(L) -2v \big) (\phi(L)+v)
\big\lfloor \frac{K - \alpha - 2}{\phi(L)+v} \big\rfloor
\Big(\big\lfloor \frac{K - \alpha - 2}{\phi(L)+v} \big\rfloor  -1 \Big)
  }{2}
\\
& \qquad \qquad +
\big(L - 2 \phi(L)-2v \big) \left\lfloor \frac{K - \alpha - 2}{\phi(L)+v} \right\rfloor \left( K-\alpha -1 - (\phi(L)+v) \left\lfloor \frac{K - \alpha - 2}{\phi(L)+v} \right\rfloor \right).~~\Box
\end{align*}

In the next section we investigate asymptotic bounds of $m$ in terms of $D$ and $N$.  This will, among other things, allow us to judge the quality of our matrices with respect to the lower bound in part 6 of Theorem~\ref{thm:MainResult}.

\subsection{Asymptotic Upper and Lower Bounds}

We begin this section by proving an asymptotic lower bound for the number of rows in any $(K,\alpha)$-coherent matrix created as per Section~\ref{sec:makeCoherent}.  Recall that we have fixed $K$ to be multiple of $\alpha$ so that $K = K_{\alpha} = \lceil D \alpha \rceil$ for some $D \in (1, \infty)$.  We have the following lower bound for $m$ as a function of $D$ and $N$.

\begin{lem}
Suppose that $2 \le D \le N^{1-\tau}$, where $\tau>0$ is some fixed constant.
For any solution to the optimization problem in Figure~\ref{fig:Problem}, where $\alpha$ can freely be chosen, one has
 $$m \gg \frac{D^2 (\log N)^2}{ \log (D \log N)}
$$
for sufficiently large values of $D \log N$.
\label{lem:AsymLowermBound}
\end{lem}

\noindent \textit{Proof:}\\

Let $Q = D \log N.$  Suppose that $S=\{s_1, \ldots s_K\}$ is a solution to the optimization problem in Figure~\ref{fig:Problem}, where $\alpha$ can be freely chosen and $K=K_{\alpha} = \lceil D \alpha \rceil$.   By Properties (I) and (II) in Figure~\ref{fig:Problem},
$$\sum_{i=1}^{K} \log s_i > \left\lfloor \frac{D \alpha}{\alpha+1} \right\rfloor  \log N \ge \frac{Q}{4}.$$
Let $q=q_S$ be the largest natural number such that $p_q \mid \prod_{i=1}^K s_i$.  For $1 \le i \le q,$ let $w_i$ be the integer such that $p_i^{w_i} \| \prod_{i=1}^K s_i$.
Since $a + b \le ab$ for $a,b \in \mathbb{N}$, it follows from Lemma \ref{lem:Sumupperbound_prime} that
$$\sum_{i=1}^q p_i^{w_i} \le \sum_{i=1}^K s_i \ll \frac{Q^2}{\log Q}$$ as $Q \to \infty$.  This implies that for $1 \le i \le q$, we have $p_i^{w_i} \le \frac{C Q^2}{\log Q}$, where $C$ is some absolute positive constant.  Therefore,
\begin{align*}
\sum_{p_i \le \sqrt{\frac{C Q^2}{\log Q}}} w_i \log p_i & \le \sum_{w=2}^\infty \sum_{p \le \left(\frac{C Q^2}{\log Q}\right)^{1/w}} w \log p \ll \sum_{2 \le w \le \frac{\log (C Q^2) -\log \log Q}{\log 2}} \frac{w Q^{2/w}}{(\log Q)^{1/w}} \ll \frac{Q}{\sqrt{\log Q}}.
\end{align*}
Since $\sum_{i=1}^q w_i \log p_i = \sum_{i=1}^K \log s_i \ge \frac{Q}{4},$ we have that
$$\sum_{p_i > \sqrt{\frac{C Q^2}{\log Q}}} w_i \log p_i \gg Q$$ for sufficiently large values of $Q$.  Let $$W = \sum_{p_i > \sqrt{\frac{C Q^2}{\log Q}}} w_i.$$  We have that $w_i \in \{0,1\}$ when $p_i > \sqrt{\frac{C Q^2}{\log Q}}$.  Also, whenever $w_i \ge 1$, it follows that $\log p_i \le \log \frac{C Q^2}{\log Q} \ll \log Q $.  Thus, for sufficiently large values of $Q$, we have $W \ge \frac{C' Q}{\log Q}$, where $C'$ is some absolute positive constant.  By the Prime Number Theorem,
$$ \sum_{i=1}^K s_i \ge \sum_{i=1}^q p_i^{w_i} \ge \sum_{i \le \frac{C' Q}{\log Q}} p_i \gg \sum_{i \le \frac{C' Q}{\log Q}} i \log i \gg \frac{Q^2}{\log Q} = \frac{D^2 (\log N)^2}{ \log (D \log N )}$$
for sufficiently large values of $Q=D \log N.~~\Box$\\

Part 6 of Theorem~\ref{thm:MainResult} informs us that $m$ must be $\Omega \left( D^2 \log_D N \right)$ for any $m \times N$ $\left(K = \lceil D \alpha \rceil,\alpha \right)$-coherent matrix.  On the other hand, Lemma~\ref{lem:AsymLowermBound} above tells that the any $m \times N$ $\left( K,\alpha \right)$-coherent matrix constructed via Section~\ref{sec:makeCoherent} must have $m = \Omega \left( \frac{D^2 (\log N)^2}{ \log (D \log N)} \right)$.  
Note that the lower bounds for matrices constructed as per Section~\ref{sec:makeCoherent} are worse by approximately a factor of $\log N$.
This is probably an indication that the $\left( K = \lceil D \alpha \rceil,\alpha \right)$-coherent matrix construction in Section~\ref{sec:makeCoherent} is suboptimal.  Certainly suboptimality of the construction in Section~\ref{sec:makeCoherent} would not be surprising given that the construction is addressing a more constrained design problem (i.e., we demand small Fourier sampling requirements).

Next we show that the asymptotically best main term for $m$ in the optimization problem in Figure~\ref{fig:Problem} can be obtained by taking each $s_j$ to be a prime.  This proves that the asymptotic lower bound given in Lemma~\ref{lem:AsymLowermBound} is tight.

\begin{lem}
\label{lem:Sumupperbound_prime}
Suppose that $2 \le D \le N^{1-\tau}$, where $\tau>0$ is some fixed constant.
If we are able to select the value for $\alpha$, the optimization problem in Figure~\ref{fig:Problem} can be solved by taking the $s_j$ to be primes in such a way that guarantees that $$m \ll \frac{D^2 (\log N)^2}{  \log (D \log N )}
$$
as $D \log N \to \infty$.

\end{lem}

\noindent \textit{Proof:}\\

Let $Q= D \log N$.
%We will consider
%asymptotic bounds as $Q \to \infty$.
Since we restrict to the case that $1 \le D \le N^{1-\tau}$, it follows that as $Q \to \infty$, we  also have that $N \to \infty$.
Let $r = \max \left( \left\lceil \frac{Q}{\log Q} \right\rceil, 9 \right).$  Note that $\log r > 2$.    Also,  by the Prime Number Theorem, $p_{r+1}  \sim Q \leq N^{1-\tau} \log N$ as $Q \to \infty$.  We will assume that $Q$ is large enough that $p_{r+1} < N$.  Choose $\alpha \in \mathbb{N}$ such that $\prod^{\alpha}_{j=1} p_{r+j} < N \leq \prod^{\alpha+1}_{j=1} p_{r+j}$.  For $1 \le i \le K_\alpha = \lceil D \alpha \rceil$, let $s_i = p_{r+i}$.  Note that our elements $s_i$ already satisfy the conditions in  Figure~\ref{fig:Problem}.  We
are left to establish a bound on $\alpha$
and then estimate $\sum_{i=1}^{K_\alpha} s_i$.

Note that $\alpha \le \beta$ whenever $\prod_{j=1}^{\beta+1} p_{r+j} \ge N$, which is equivalent to $\sum_{j=1}^{\beta+1} \log p_{r+j} \ge \log N$.  Let $\beta = \left\lceil \frac{2 \log N}{\log r} \right\rceil.$
We have that $p_k \ge k$ for $k\ge 1$.
Hence, we have that
\begin{align*}
\sum_{i=1}^{\beta+1} \log p_{r+i}
& \ge \sum_{i=1}^{\beta+1} \log (r+i)
\\ & \ge \int_{r}^{r+\beta+1} \log x \, dx
\\ & = (r+\beta + 1) \log (r + \beta + 1) - (r + \beta + 1) - r \log r + r
\\ & = (\beta + 1) (\log (r + \beta  + 1) - 1) +r \log \big(1 + \frac{\beta +1}{r}\big)
\\ & > (\beta + 1) (\log (r + \beta  + 1) - 1)
\\ & > \beta (\log r -1)
\\ & > \frac{2 \log N}{\log r} \cdot \frac{\log r}{2}
\\ & = \log N.
\end{align*}
Note that as $Q \to \infty$, $r+K_\beta \ll \frac{Q}{\log Q}$.
Thus, since $\alpha \le \beta$, we have by \cite[Lemma 6]{CISS} and the Prime Number Theorem that
\begin{equation*}
\sum_{i=1}^{K_\alpha} p_{r+i}  \le \sum_{i=1}^{K_\beta} p_{r+i}
  \le C_1 \frac{p_{r+K_\beta}^2}{\log p_{r+K_\beta}}
\le C_2 \frac{\big((r+K_\beta) \log(r+K_\beta) \big)^2}{\log \big( (r+K_\beta) \log(r+K_\beta) \big)}
,
\end{equation*}
for some absolute constants $C_1$ and $C_2$. As $Q \to \infty$,
\begin{align*}
m \ll \frac{\big((r+K_\beta) \log(r+K_\beta) \big)^2}{\log \big( (r+K_\beta) \log(r+K_\beta) \big)}
 & \ll  \frac{Q^2}{\log Q}
 =  \frac{D^2 (\log N)^2}{ \log (D \log N )}.~~\Box
\end{align*}

Although Lemma~\ref{lem:Sumupperbound_prime} shows that simply using primes for our $s_j$ values is asymptotically optimal, it is important to note that the convergence of such primes-only solutions to the optimal value as $D \log N \to \infty$ is likely very slow.  For real world values of $N$ and $D$ the more general criteria that the $s_j$ values be pairwise relatively prime can produce significantly smaller $m$ values.  This is demonstrated empirically in Section~\ref{sec:Experiments}.  However, Lemma~\ref{lem:Sumupperbound_prime} also formally justifies the idea that the $s_j$ values can be restricted to smaller subsets of relatively prime integers (e.g., the prime numbers) before solving the optimization problem in Figure~\ref{fig:Problem} without changing the asymptotic performance of the generated solutions.  This idea can help make the (approximate) solution of the optimization problem in Figure~\ref{fig:Problem} more computationally tractable in practice.

\section{Formulation of the Matrix Design Problem as a Linear Integer Program}
\label{sec:LinIntProg}

To formulate the problem as a linear integer program, we define $K = K_\alpha = \left \lceil D\alpha \right \rceil$ and $B = p_{t+K-\alpha-1}$ as in Part 4 of Lemma  \ref{lem:sjBounds}.  Let $s_{j,i} \in \{ 0,1 \}$ for $j \in [1,K] \cap \mathbbm{N}$ and $i \in [1,B] \cap \mathbb{N}$.  We then let
$\displaystyle s_j = \sum^{B}_{i=1} s_{j,i} \cdot i$ and, for $k \in [1, t+K-\alpha -1] \cap \mathbb{N}$, define
\begin{align}
\delta_{k,i} & = \left\{ \begin{array}{ll} 1 & \textrm{if } p_k \mid i \\ 0 & {\rm otherwise} \end{array} \right. .
\label{eqn:Def_LinProgr}
\end{align}
Then, for a given $\alpha$, we can minimize Equation~(\ref{eqn:Optsamps}) by minimizing
\begin{equation}
m = \sum^{K}_{j = 1} s_j = \sum^{K}_{j = 1} \sum^{B}_{i=1}  s_{j,i}  \cdot i
\label{eqn:Optsamps2}
\end{equation}
subject to the following linear constraints:
\begin{enumerate}
\item $\sum^{B}_{i=1} s_{j,i}=1$ for all $j \in [1,K] \cap \mathbbm{N}$.
\item $s_{j,i} \in \{0,1\}$ for all $j \in [1,K] \cap \mathbbm{N}$ and $i \in [1,B] \cap \mathbbm{N}$.
\item $\sum^{B}_{i=1}(i\cdot s_{j+1,i}-i \cdot s_{j,i})\geq 1$ for all $j \in [1,K-1] \cap \mathbbm{N}$.
\item $\sum^{\alpha}_{j=1} \sum^{B}_{i=1}  s_{j,i} \cdot \ln i < \ln N \leq \sum^{\alpha+1}_{j=1} \sum^{B}_{i=1} s_{j,i} \cdot \ln i$.
\item $s_{j,i}=0$ for all $j \in [1,K] \cap \mathbbm{N}$ and $i \in [1,p_i-1] \cap \mathbbm{N}$.
\item $\sum^{K}_{j = 1}\sum^{B}_{i=1} \delta_{k,i} \cdot s_{j,i} \leq 1$ for all $k \in [1, t+K-\alpha -1] \cap \mathbb{N}$.
\end{enumerate}

The first and second constraint together state that for each $j$, $s_{j,i}$ is non-zero for exactly one value of $i \in [1,B] \cap \mathbbm{N}$, implying that $\displaystyle s_j = i$.  This in turn, by the third constraint, implies that $s_1 < s_2 < \cdots < s_K$, which is Constraint I in Figure \ref{fig:Problem}.  Upon applying the natural logarithm in Constraint II in Figure \ref{fig:Problem} to convert a nonlinear constraint to a linear constraint, one obtains something equivalent to our fourth constraint above.  The fifth constraint above simply forces $s_j \geq p_j$, which will be true for any solution to the optimization problem in Figure \ref{fig:Problem}.  The last constraint ensures that $s_1, \ldots, s_K$ are pairwise relatively prime, which is Constraint III in Figure \ref{fig:Problem}.  Hence, the optimization problem in this section is equivalent to the optimization problem in Figure \ref{fig:Problem}.

\section{Numerical Experiments}
\label{sec:Experiments}

In this section we investigate the optimal Fourier sampling requirements related to $m \times N$ $\left(K = \left \lceil D\alpha \right \rceil,\alpha \right)$-coherent matrices, optimized over the $\alpha$ parameter, for several values of $D$ and $N$.  This is done for given values of $D \in (1,\infty)$ and $N \in \mathbbm{N}$ by solving the optimization problem in Figure~\ref{fig:Problem} via the linear integer program presented in Section~\ref{sec:LinIntProg} for all feasible values of $\alpha \in \left[1, \log_2 N \right] \cap \mathbbm{N}$.\footnote{It is important to note that many values of $\alpha$ can be disqualified as optimal without solving a linear integer program by comparing previous solutions to the lower bounds given in Lemma~\ref{lem:mLowerBound} and Corollary~\ref{cor:mLowerBound2}.}  The solution yielding the smallest Fourier sampling requirement, $m - K_\alpha + 1$ from Lemma~\ref{lem:SmallFourier}, for the given $D$ and $N$ values (minimized over all $\alpha$ values) is the one reported for experiments in this section.  Each linear integer program was solved with IBM ILOG OPL-CPLEX with parameters generated using Microsoft Visual Studio.  Examples of the actual files ran can be downloaded from the contact author's website.\footnote{http://www.math.duke.edu/$\sim$markiwen/DukePage/code.htm}

In order to make our numerical experiments more meaningful we computed optimal incoherent matrices which also have the RIP$_2$ (see part 3 of Theorem~\ref{thm:MainResult}).  Hence, we set $D = \frac{k-1}{\epsilon}$ for a given sparsity value $k \in [1,N] \cap \mathbbm{N}$ and $\epsilon \in (0,1)$.  In all experiments the value of $\epsilon$ was fixed to be slightly less than $3 / \left( 4 + \sqrt{6} \right) \approx 0.465$ which ensures that $l^1$-minimization can be utilized with the produced RIP$_2$ matrices for accurate Fourier approximation (e.g., see Theorem 2.7 in \cite{HolgerStableRecov}).  

\begin{figure}
\resizebox{0.5\columnwidth}{!}{\includegraphics{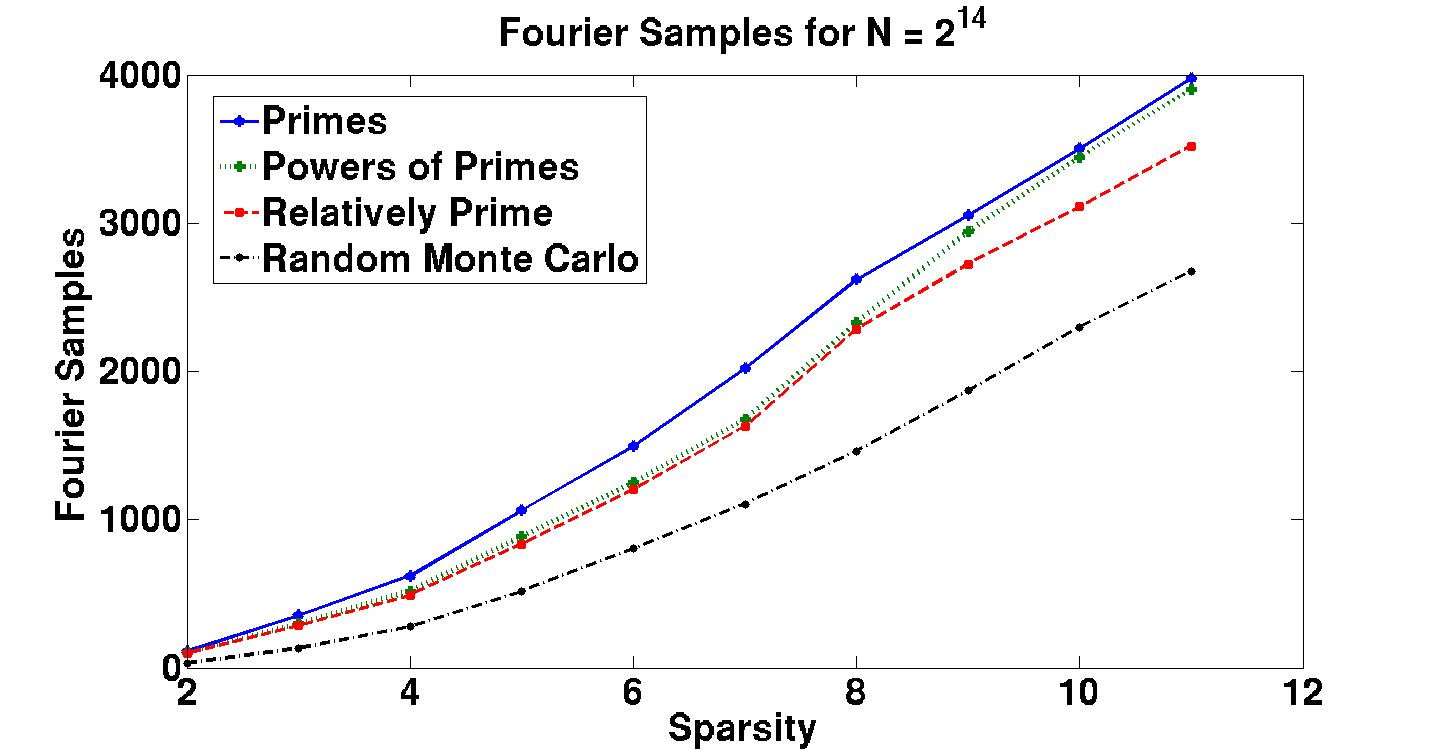}}
\resizebox{0.5\columnwidth}{!}{\includegraphics{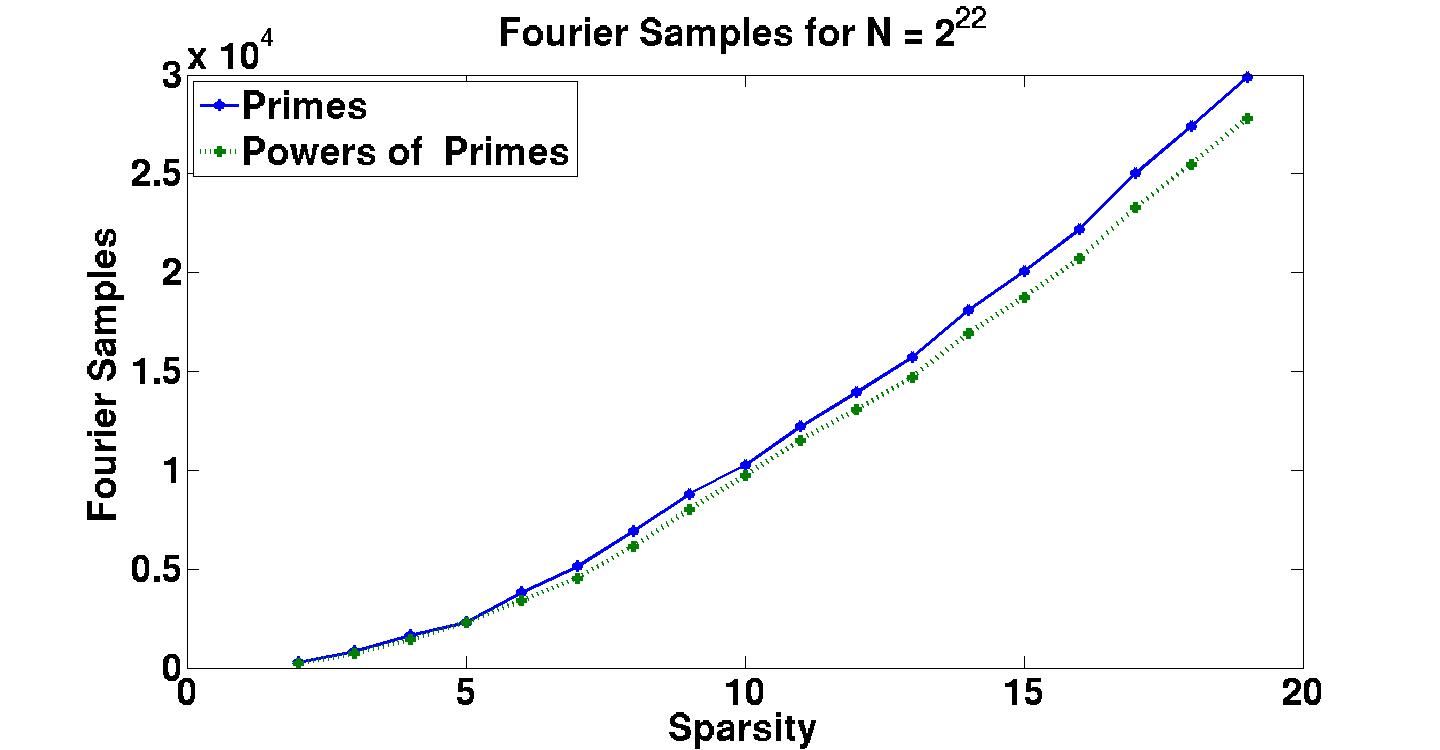}}
\caption{
\textbf{On the Left:}  The minimal Fourier sampling requirements, $m - K_\alpha + 1$ minimized over all feasible $\alpha \in [1,14] \cap \mathbbm{N}$, for any possible $m \times 2^{14}$ $\left(K = \left \lceil \frac{(k-1)\alpha}{\epsilon} \right \rceil,\alpha \right)$-coherent matrix constructed via Section~\ref{sec:makeCoherent}.  Here $\epsilon$ was fixed to be $4 / \left( 6 + \sqrt{7} \right) \approx 0.463$, and the sparsity parameter, $k$, was varied between $2$ and $11$.  \textbf{On the Right:}  The minimal Fourier sampling requirements, $m - K_\alpha + 1$ minimized over all feasible $\alpha \in [1,22] \cap \mathbbm{N}$, for any possible $m \times 2^{22}$ $\left(K = \left \lceil \frac{(k-1)\alpha}{\epsilon} \right \rceil,\alpha \right)$-coherent matrix constructed via Section~\ref{sec:makeCoherent}.  Here $\epsilon$ was again fixed to be $4 / \left( 6 + \sqrt{7} \right) \approx 0.463$, and the sparsity parameter, $k$, was varied between $2$ and $19$.
}
\label{fig:Samps}
\end{figure}

Three variants of the optimization problem in Figure~\ref{fig:Problem} were solved in order to determine the minimal Fourier sampling requirements associated with various classes of $\left(\left \lceil \frac{(k-1)\alpha}{\epsilon} \right \rceil,\alpha \right)$-coherent matrices created via Section~\ref{sec:makeCoherent}.  These three variants include the:  
\begin{enumerate}
\item \textbf{Relatively Prime} optimization problem exactly as stated in Figure~\ref{fig:Problem} and reformulated in Section~\ref{sec:LinIntProg}.
\item \textbf{Powers of Primes} optimization problem.  Here the $s_j$ values are further restricted to each be a power of a single prime number.  
\item \textbf{Primes} optimization problem.  Here each $s_j$ value is further restricted to simply be a prime number.
\end{enumerate}
These different variants allow some trade off between computational complexity and the minimality of the generated incoherent matrices.  See Figure~\ref{fig:Samps} for a comparison of the solutions to these optimization problems for two example values of $N$.

In creating the solutions graphed in Figure~\ref{fig:Samps} computer memory was the primary constraining factor.  For each of the two values of $N$ the sparsity, $k$, was increased until computer memory began to run out during the solution of one of the required linear integer programs.\footnote{A modest desktop computer with an Intel Core i7-920 processor @ 2.67 Ghz and 2.99 GB of RAM was used to solve all linear integer programs reported on in Figure~\ref{fig:Samps}.}  All linear integer programs which ran to completion did so in less than 90 minutes (most finishing in a few minutes or less).  Not surprisingly, the relatively prime solutions always produce smaller Fourier sampling requirements than the more restricted powers of primes solutions, with the tradeoff being that they are generally more difficult to solve.  Similarly, the powers of primes solutions always led to smaller Fourier sampling requirements than the even more restricted primes solutions.

For the sake of comparison, the left plot in Figure~\ref{fig:Samps} also includes Fourier sampling results for RIP$_2$($2^{14}$,$k$,$\epsilon<0.465$) matrices created via random sampling based incoherence arguments for each sparsity value.  These random Fourier sampling requirements were calculated by choosing rows from an $2^{14} \times 2^{14}$ inverse DFT matrix, $\mathcal{F}^{-1}$, uniformly at random without replacement.  After each row was selected, the $\mu$-coherence of the submatrix formed by the currently selected rows was calculated (see Definition~\ref{def:Coherence}).  As soon as the coherence became small enough that Theorem~\ref{thm:RIPcoherent} guaranteed that the matrix would have the RIP$_2$($2^{14}$,$k$,$\epsilon<0.465$) for the given value of $k$, the total number of inverse DFT rows selected up to that point was recored as a trial Fourier sampling value.  This entire process was repeated 100 times for each value of $k$.  The smallest Fourier sampling value achieved out of these 100 trials was then reported for each sparsity $k$ in the left plot of Figure~\ref{fig:Samps}.

Looking at the plot of the left in Figure~\ref{fig:Samps} we can see that the randomly selected submatrices guaranteed to have the RIP$_2$ require fewer Fourier samples than the deterministic matrices constructed herein.  Hence, if Fourier sampling complexity is one's primary concern, traditional matrix design techniques should be utilized.  However, it is important to note that such randomly constructed Fourier matrices cannot currently be utilized in combination with $o(N)$-time Fourier approximation algorithms.  Our deterministic incoherent matrices, on the other hand, have associated sublinear-time approximation algorithms (see the first part of Theorem~\ref{thm:MainResult}).

Finally, we conclude this paper by noting that heuristic solutions methods can almost certainly be developed for solving the optimization problem in Figure~\ref{fig:Problem}.  Such methods are often successful at decreasing memory usage and computation time while still producing near-optimal results.  We leave further consideration of such approaches as future work.

\bibliographystyle{abbrv}
\bibliography{ripmeRefs}

\appendix

\section{Proof of the First Statement in Theorem~\ref{thm:MainResult}}
\label{appsec:SublinearProof}

We will prove a slightly more general variant of the first statement given in Theorem~\ref{thm:MainResult}.  Below we will work with $\left( K, c_{\rm min}, \alpha \right)$-coherent matrices.

\begin{Definition}
Let $K \in [1,N] \cap \mathbbm{N}$ and $c_{\rm min}, \alpha \in \mathbbm{R}^{+}$.  An $m \times N$ positive real matrix, $\mathcal{M} \in [0,\infty)^{m \times N}$, is called $\left( K, c_{\rm min}, \alpha \right)$-coherent if both of the following properties hold:
\begin{enumerate}
\item Every column of $\mathcal{M}$ contains at least $K$ nonzero entries. 
\item All nonzero entries are at least as large as $c_{\rm min}$.
\item For all $j, l \in [0,N)$ with $j \neq l$, the associated columns, $\mathcal{M}_{\cdot,j} ~{\rm and}~ \mathcal{M}_{\cdot,l} \in [0,\infty)^m$, have $\mathcal{M}_{\cdot,j} \cdot \mathcal{M}_{\cdot,l} \leq \alpha$.
\end{enumerate}
\label{appdef:PosCoherent}
\end{Definition}

Clearly, any $\left( K, \alpha \right)$-coherent matrix will also be $\left( K, 1, \alpha \right)$-coherent.  Other examples of $\left( K, c_{\rm min}, \alpha \right)$-coherent matrices include ``corrupted'' or ``noisy'' $\left( K, \alpha \right)$-coherent matrices, as well as matrices whose columns are spherical code words from the first orthant of $\mathbbm{R}^m$.  In what follows we will generalize results and constructions from \cite{ImpFourier}.  We will give self-contained proofs whenever possible, although it will be necessary on occasion to state generalized results from \cite{ImpFourier} whose proofs we omit.

\subsection{Some Useful Properties of $\left( K, c_{\rm min}, \alpha \right)$-Coherent Matrices}
\label{appsec:MProps}

In what follows, $\mathcal{M} \in [0,\infty)^{m \times N}$ will always refer to a given $m \times N$ $\left( K, c_{\rm min}, \alpha \right)$-coherent matrix.  Let $n \in [0,N) \cap \mathbbm{N}$.  We define $\mathcal{M}(K,n)$ to be the $K \times N$ matrix created by selecting the $K$ rows of $\mathcal{M}$ with the largest entries in its $n^{\rm th}$ column.  Furthermore, we define $\mathcal{M'}(K,n)$ to be the $K \times (N-1)$ matrix created by deleting the $n^{\rm th}$ column of $\mathcal{M}(K,n)$.  Thus, if
$$\mathcal{M}_{j_1,n} ~\geq~ \mathcal{M}_{j_2,n} ~\geq~ \dots ~\geq~ \mathcal{M}_{j_m,n}$$
then
\begin{equation}
\mathcal{M}(K,n) = \left( \begin{array}{l} \mathcal{M}_{j_1} \\ \mathcal{M}_{j_2} \\ \vdots \\ \mathcal{M}_{j_K} \\ 
\end{array} \right)
\label{appdef:Msub}
\end{equation}
and
\begin{equation}
\mathcal{M'}(K,n) = \left( \begin{array}{lllllll} \mathcal{M}_{j_1,1} & \mathcal{M}_{j_1,2} & \dots & \mathcal{M}_{j_1,n-1} & \mathcal{M}_{j_1,n+1} & \dots & \mathcal{M}_{j_1,N} \\ 
\mathcal{M}_{j_2,1} & \mathcal{M}_{j_2,2} & \dots & \mathcal{M}_{j_2,n-1} & \mathcal{M}_{j_2,n+1} & \dots & \mathcal{M}_{j_2,N} \\ 
&&&&\hspace{-5pt} \vdots&& \\ 
\mathcal{M}_{j_K,1} & \mathcal{M}_{j_K,2} & \dots & \mathcal{M}_{j_K,n-1} & \mathcal{M}_{j_K,n+1} & \dots & \mathcal{M}_{j_K,N} \\ 
\end{array} \right).
\label{appdef:M'sub}
\end{equation}
The following two lemmas motivate the main results of this section.

\begin{lem}
Suppose $\mathcal{M}$ is a $\left(K, c_{\rm min} ,\alpha \right)$-coherent matrix.  Let $n \in [0,N) \cap \mathbbm{N}$, $k \in \left[1,K \cdot c_{\rm min}^2 / \alpha \right] \cap \mathbbm{N}$, and $\vec{x} \in \mathbbm{C}^{N-1}$.  Then, at most $\frac{k \alpha}{c_{\rm min}^2}$ of the $K$ entries of $\mathcal{M'}(K,n) \cdot \vec{x}$ will have magnitude greater than or equal to $\frac{c_{\rm min}}{k} \cdot \| \vec{x} \|_1$.
\label{applem1}
\end{lem}

\noindent \textit{Proof:}\\  

We have that
$$\left| \left\{ j ~\bigg|~ \left| \left(\mathcal{M'}(K,n) \cdot \vec{x} \right)_j \right| \geq \frac{c_{\rm min} \cdot \| \vec{x} \|_1 }{k} \right\} \right| \leq \frac{k \cdot \left\| \mathcal{M'}(K,n) \cdot \vec{x} \right\|_1 }{c_{\rm min} \cdot \| \vec{x} \|_1} \leq \frac{k}{c_{\rm min}} \cdot \| \mathcal{M'}(K,n) \|_1.$$  
Focusing now on $\mathcal{M'}(K,n)$, we can see that
\begin{equation}
\left \| \mathcal{M'}(K,n) \right\|_1 = \max_{l \in [1,N-1] \cap \mathbbm{N}} \left\| \left(\mathcal{M'}(K,n) \right)_{\cdot,l} \right\|_1 \leq \max_{l \in [1,N-1] \cap \mathbbm{N}} \frac{\left\langle \left(\mathcal{M}(K,n) \right)_{\cdot,n}, \left(\mathcal{M}(K,n) \right)_{\cdot,l} \right\rangle}{c_{\rm min}} \leq \frac{\alpha}{c_{\rm min}}.
\label{appeqn:l1Bound}
\end{equation}
The result follows.~~$\Box$ \\

With Lemma~\ref{applem1} in hand we can now prove our second lemma.  However, we must first establish some notation.  For any given $\vec{x} \in \mathbbm{C}^{N}$ and subset $S \subseteq [0,N) \cap \mathbbm{N}$, we will let $\vec{x}_{S} \in \mathbbm{C}^{N}$ be equal to $\vec{x}$ on the indexes in $S$ and be zero elsewhere.  Thus, 
$$\left( \vec{x}_{S} \right)_i = \left\{ \begin{array}{ll} x_i & \textrm{if } i \in S \\ ~0 & {\rm otherwise} \end{array} \right..$$
Furthermore, for a given integer $k < N$, we will let $S^{\rm opt}_k \subset [0,N) \cap \mathbbm{N}$ be the first $k$ element subset of $[0,N) \cap \mathbbm{N}$ in lexicographical order with the property that $|x_s| \geq |x_t|$ for all $s \in S^{\rm opt}_k$ and $t \in [0,N) \cap \mathbbm{N} - S^{\rm opt}_k$.  Thus, $S^{\rm opt}_k$ contains the indexes of $k$ of the largest magnitude entries in $\vec{x}$.  Finally, we will define $\vec{x}^{\rm opt}_k$ to be $\vec{x}_{S^{\rm opt}_k}$, an optimal $k$-term approximation to $\vec{x}$.

\begin{lem}
Suppose $\mathcal{M}$ is a $\left(K, c_{\rm min} ,\alpha \right)$-coherent matrix.  Let $n \in [0,N) \cap \mathbbm{N}$, $\tilde{k} \in \left[1,K \cdot c^2_{\rm min} / \alpha \right] \cap \mathbbm{N}$, $S \subset [0,N) \cap \mathbbm{N}$ with $|S| \leq \tilde{k}$, and $\vec{x} \in \mathbbm{C}^{N-1}$.  Then, $\mathcal{M'}(K,n) \cdot \vec{x}$ and $\mathcal{M'}(K,n) \cdot \left( \vec{x} - \vec{x}_{S} \right)$ will differ in at most $\frac{\tilde{k} \alpha}{c_{\rm min}^2}$ of their $K$ entries.  
\label{applem2}
\end{lem}

\noindent \textit{Proof:} \\

Let $\vec{\mathbbm{1}} \in  \mathbbm{C}^{N-1}$ be the vector of all ones.  We have that
$$\left| \left\{ j ~\bigg|~ \left(\mathcal{M'}(K,n) \cdot \vec{x} \right)_j  \neq \left( \mathcal{M'}(K,n) \cdot \left( \vec{x} - \vec{x}_{S} \right) \right)_j \right\} \right| = \left| \left\{ j ~\bigg|~ \left(\mathcal{M'}(K,n) \cdot \vec{x}_{S} \right)_j  \neq 0 \right\} \right| \leq \left| \left\{ j ~\bigg|~ \left(\mathcal{M'}(K,n) \cdot \vec{\mathbbm{1}}_{S} \right)_j  \geq c_{\rm min} \right\} \right|$$
since all the nonzero entries of $\mathcal{M'}(K,n)$ are at least as large as $c_{\rm min}$.  Applying Lemma~\ref{applem1} with $\vec{x} = \vec{\mathbbm{1}}_{S}$ and $k = \left \| \vec{\mathbbm{1}}_{S} \right \|_1 = |S|$ finishes the proof.~~$\Box$ \\

By combining the two Lemmas above we are able to bound the accuracy with which we can approximate any entry of an arbitrary complex vector $\vec{x} \in \mathbbm{C}^N$ using only the measurements from a $\left(K, c_{\rm min} ,\alpha \right)$-coherent matrix.  The following lemma motivates the remainder of this appendix.

\begin{lem}
Suppose $\mathcal{M}$ is a $\left(K, c_{\rm min} ,\alpha \right)$-coherent matrix.  Let $n \in [0,N) \cap \mathbbm{N}$, $k \in \left[1,K \cdot c^2_{\rm min} / \alpha \right] \cap \mathbbm{N}$, $\epsilon \in (0,1]$, $c \in [2,\infty) \cap \mathbbm{N}$, and $\vec{x} \in \mathbbm{C}^{N}$.  If $K > c \cdot (k \alpha/ c_{\rm min}^2 \epsilon)$ then more than $\frac{c-2}{c} \cdot K$ of the $K$ entries of $\mathcal{M}(K,n) \cdot \vec{x}$ can be used to estimate $x_n$ to within $\frac{\epsilon \left\| \vec{x} - \vec{x}^{\rm opt}_{(k/\epsilon)} \right\|_1}{k}$ accuracy.
\label{applem:Mest}
\end{lem}

\noindent \textit{Proof:}\\  

Define $\vec{y} \in \mathbbm{C}^{N-1}$ to be $\vec{y} = \left( x_0, x_1, \dots, x_{n-1}, x_{n+1}, \dots, x_{N-1} \right)$.  We have
$$\mathcal{M}(K,n) \cdot \vec{x} ~=~ x_n \cdot \left(\mathcal{M}(K,n)\right)_{\cdot,n} + \mathcal{M'}(K,n) \cdot \vec{y}.$$  Applying Lemma~\ref{applem2} with $\tilde{k} = (k/\epsilon)$ demonstrates that at most 
$\frac{k \alpha}{\epsilon \cdot c_{\rm min}^2}$ entries of $\mathcal{M'}(K,n) \cdot \vec{y}$ differ from their corresponding entries in $\mathcal{M'}(K,n) \cdot \left( \vec{y} - \vec{y}^{\rm opt}_{(k / \epsilon)} \right)$.  
Of the remaining $K - \frac{k \alpha}{\epsilon \cdot c_{\rm min}^2}$ entries of $\mathcal{M'}(K,n) \cdot \vec{y}$, at most $\frac{k \alpha}{\epsilon \cdot c_{\rm min}^2}$ will have magnitudes greater than or equal to 
$\epsilon \cdot c_{\rm min} \left\| \vec{y} - \vec{y}^{\rm opt}_{(k / \epsilon)} \right\|_1 / k$ by Lemma~\ref{applem1}.  Hence, at least 
$$K - 2 \left( \frac{k \alpha}{\epsilon \cdot c_{\rm min}^2} \right) > \frac{c-2}{c} \cdot K$$
entries of $\mathcal{M'}(K,n) \cdot \vec{y}$ will have a magnitude no greater than
$$\frac{\epsilon \cdot c_{\rm min} \left\| \vec{y} - \vec{y}^{\rm opt}_{(k / \epsilon)} \right\|_1}{k} \leq \frac{\epsilon \cdot c_{\rm min} \left\| \vec{x} - \vec{x}^{\rm opt}_{(k / \epsilon)} \right\|_1}{k}.$$
Therefore, $\frac{\left(\mathcal{M}(K,n) \cdot \vec{x}\right)_j}{\left(\mathcal{M}(K,n)\right)_{j,n}}$ will approximate $x_n$ to within the stated accurracy for more than $\frac{c-2}{c} \cdot K$ values $j \in [1,K] \cap \mathbbm{N}$.~~$\Box$ \\

Lemma~\ref{applem:Mest} generalizes Theorem 4 in Section 3 of \cite{ImpFourier}.  Thus, Lemma~\ref{applem:Mest} can be used to modify the proof of Theorem 6 in Section 4 of \cite{ImpFourier} in order to prove that a variant of Algorithm 2 from \cite{ImpFourier} will provide instance optimal approximation guarantees along the lines of Equation~\ref{eqn:Aerror2} for general compressed sensing recovery problems.\footnote{More precisely, the error bound in Equation 21 of Theorem 6 in \cite{ImpFourier} holds without the additional third $22 \sqrt{k} ~ \left\| \cdot \right\|_1$ term.}  The intuitive idea is as follows:  If the constant $c$ from Lemma~\ref{applem:Mest} is set to be at least $4$, then more than half of the $K$ entries of $\mathcal{M}(K,n) \cdot \vec{x}$ can accurately estimate $x_n$.  This is enough to guarantee that the imaginary part of $x_n$ will be accurately estimated by the median of the imaginary parts of all $K$ properly scaled entries of $\mathcal{M}(K,n) \cdot \vec{x}$.  Of course, the real part of $x_n$ can also be estimated in a similar fashion.  Hence, given both $\mathcal{M}$ and $\mathcal{M} \vec{x} \in \mathbbm{C}^m$, Lemma~\ref{applem:Mest} ensures that computing $N$ medians of $K$ reweighted elements of $\mathcal{M} \vec{x}$ will allow us to accurately estimate all $N$ entries of $\vec{x}$.  If we do this and then report only the largest $2k$ estimates in magnitude, together with their vector indexes, we will obtain an approximation for $\vec{x}$ which is at least as good as $\vec{x}_{S^{\rm opt}_k} ~=~ \vec{x}^{\rm opt}_k$.  See Algorithm 2 and Theorem 4 in \cite{ImpFourier} for a detailed proof in the Fourier setting.

It is worth noting that the randomized approximation results in \cite{ImpFourier} also generalize in this manner (i.e., Corollaries 3 and 4 in \cite{ImpFourier}).  If a small set of rows is randomly selected from a $\left(K, c_{\rm min} ,\alpha \right)$-coherent matrix, the resulting submatrix can still be used to yield an accurate instance optimal approximation for any $\vec{x} \in \mathbbm{C}^N$ with high probability.  The proof of this fact follows from Corollary~\ref{appcor:RandM} below.  However, before we can state the corollary we need an additional definition:  For any multiset, $$\tilde{s} = \left\{ \tilde{s}_1,~\tilde{s}_2,\dots,\tilde{s}_\beta \right\} \subset [1,m] \cap \mathbbm{N},$$  
we will let $\mathcal{M}_{\tilde{s}}$ denote the $\beta \times N$ matrix formed by the $\beta$ rows of $\mathcal{M}$ listed in $\tilde{s}$.  In other words, 
\begin{equation}
\mathcal{M}_{\tilde{s}} = \left( \begin{array}{l} \mathcal{M}_{\tilde{s}_1} \\ \mathcal{M}_{\tilde{s}_2} \\ \vdots \\ \mathcal{M}_{\tilde{s}_\beta} \\ 
\end{array} \right).
\label{appdef:RandMsub}
\end{equation} 
We have the following result.

\begin{cor}
Suppose $\mathcal{M}$ is an $m \times N$ $\left(K, c_{\rm min} ,\alpha \right)$-coherent matrix.  Let $k \in \left[1,K \cdot c^2_{\rm min} / \alpha \right] \cap \mathbbm{N}$, $\epsilon \in (0,1]$, $c \in [14,\infty) \cap \mathbbm{N}$, $\sigma \in [2/3,1)$, and $\vec{x} \in \mathbbm{C}^{N}$.  Select a multiset of the rows of $\mathcal{M}$, $\tilde{s} \subset [1,m] \cap \mathbbm{N}$, by independently choosing 
\begin{equation}
\beta \geq 28.56 \cdot \frac{m}{K} \ln \left( \frac{2N}{1-\sigma} \right)
\label{appeqn:Rand_rows}
\end{equation}
values from $[1,m] \cap \mathbbm{N}$ uniformly at random with replacement.  If $K > c \cdot (k \alpha/ c_{\rm min}^2 \epsilon)$ then $\mathcal{M}_{\tilde{s}}$ will have both of the following properties with probability at least $\sigma$:  
\begin{enumerate}
\item There will be at least $\tilde{l} = 21 \cdot \ln \left( \frac{2N}{1-\sigma} \right)$ nonzero values in every column of $\mathcal{M}_{\tilde{s}}$.  Hence, $\mathcal{M}_{\tilde{s}}(\tilde{l},n)$ will be well defined for all $n \in [0,N) \cap \mathbbm{N}$ (see Equation~\ref{appdef:Msub}).
\item For all $n \in [0,N) \cap \mathbbm{N}$ more than $\tilde{l} / 2$ of the entries in $\mathcal{M}_{\tilde{s}}(\tilde{l},n) \cdot \vec{x}$ ~(i.e., more than half of the values $j \in [1,\tilde{l}] \cap \mathbbm{N}$, counted with multiplicity) will have
$$\left| \frac{\left(\mathcal{M}_{\tilde{s}}(l,n) \cdot \vec{x} \right)_j}{\left(\mathcal{M}_{\tilde{s}}(l,n) \right)_{j,n}} - x_n \right| \leq \frac{\epsilon \left\| \vec{x} - \vec{x}^{\rm opt}_{(k/\epsilon)} \right\|_1}{k}.$$
\end{enumerate}
\label{appcor:RandM}
\end{cor}

\noindent \textit{Proof:}\\  

Fix $n \in [0,N) \cap \mathbbm{N}$.  We select our multiset, $\tilde{s} \subset [1,m] \cap \mathbbm{N}$, of the rows of $\mathcal{M}$ by choosing $\beta$ elements of $[1,m] \cap \mathbbm{N}$ uniformly at random with replacement.  
Denote the $j^{\rm th}$ element chosen for $\tilde{s}$ by $\tilde{s}_j$.  Finally, let $P^n_j$ be the random variable indicating whether $\mathcal{M}_{\tilde{s}_j,n} > 0$, and let $Q^n_j$ be the random variable indicating whether $\tilde{s}_{j}$ satisfies
\begin{equation}
\left| \frac{\left(\mathcal{M} \cdot \vec{x}\right)_{\tilde{s}_j}}{\mathcal{M}_{\tilde{s}_j,n}} - x_n \right| \leq \frac{\epsilon \left\| \vec{x} - \vec{x}^{\rm opt}_{(k/\epsilon)} \right\|_1}{k}
\label{appeqn:PropRand}
\end{equation}
conditioned on $P^n_j$.  Thus, $P^n_j = 1$ if $\mathcal{M}_{\tilde{s}_j,n} > 0$, and $0$ otherwise.  Similarly, 
$$Q^n_j = \left\{ \begin{array}{ll} 1 & \textrm{if } \tilde{s}_{j} \textrm{ satisfies Equation~\ref{appeqn:PropRand} and } P^n_j = 1\\ 0 & {\rm otherwise} \end{array} \right..$$
Lemma~\ref{applem:Mest} implies that $\mathbbm{P}\left[ Q^n_j = 1 ~\big|~ P^n_j = 1 \right] > \frac{6}{7}$.  Furthermore, $\mu = \mathbbm{E} \left[ \sum^{\beta}_{j = 1} Q^n_j ~\big|~ P^n_1, \dots, P^n_{\beta} \right] \geq \frac{6}{7} \left( \sum^{\beta}_{j = 1} P^n_j \right)$.

Let $l = \sum^{\beta}_{j = 1} P^n_j$.  The Chernoff bound guarantees that 
$$\mathbbm{P}\left[\sum^{\beta}_{j = 1} Q^n_j < \frac{4 \cdot l}{7} ~\bigg|~ l \right] \leq e^{- \frac{\mu}{18}} \leq e^{-\frac{l}{21}}.$$  
Thus, if $l > 21$ we can see that $\sum^{\beta}_{j = 1} Q^n_j$ will be less than $\frac{l + 1}{2}$ with probability less than $e^{-\frac{l}{21}}$.  Hence, if $l \geq 21 \ln \left( \frac{2N}{1-\sigma} \right),$ then Property~2 will fail to be satisfied for $n$ with probability less than $\frac{1-\sigma}{2N}$.  Focusing now on $l$, we note that $\mathbbm{P}\left[ P^n_j = 1 \right] \geq \frac{K}{m}$ so that $\tilde{\mu} = \mathbbm{E} \left[ l \right] \geq \frac{K}{m} 
\beta$.

Let $\tilde{l} = 21 \ln \left( \frac{2N}{1-\sigma} \right)$. Applying the Chernoff bound one additional time reveals that $\mathbbm{P}\left[ l ~<~ \tilde{l} \right] < \mathbbm{e}^{- \tilde{\mu} \cdot \left( 1 - \frac{\tilde{l}}{\tilde{\mu}}  \right)^2 / 2 }.$  Hence, if we wish to bound $\mathbbm{P}\left[ l ~<~ \tilde{l} \right]$ from above by $\frac{1-\sigma}{2N}$ it suffices to have $\tilde{\mu}^2 - \frac{44}{21} \tilde{\mu} \tilde{l} + \tilde{l}^2 \geq 0$.  Setting $\beta \geq 1.36 \cdot \frac{m}{K} \tilde{l} = 28.56 \cdot \frac{m}{K} \ln \left( \frac{2N}{1-\sigma} \right)$ achieves this goal.  The end result is that $\mathcal{M}_{\tilde{s}}$ will fail to satisfy both Properties~1 and~2 for any $n \in [0,N) \cap \mathbbm{N}$ with probability less than $\frac{1-\sigma}{N}$.  Applying the union bound over all $n \in [0,N) \cap \mathbbm{N}$ finishes the proof.~~$\Box$ \\

Corollary~\ref{appcor:RandM} considers selecting a multiset of rows from a $\left(K, c_{\rm min} ,\alpha \right)$-coherent matrix.  Hence, some rows may be selected more than once.  If this occurs, rows should be considered to be selected multiple times for counting purposes only.  That is, all computations involving a row which is selected several times should still be carried out only once.  However, the results of these computations should be considered with greater weight during subsequent reconstruction efforts (e.g., multiply selected rows should be considered as generating multiple duplicate entries in $\mathcal{M}_{\tilde{s}} \cdot \vec{x}$).  

In this paper we are primarily concerned with guaranteed approximation results.  Hence, we will leave further consideration of randomized approximation techniques to the reader.  Instead, we will now consider fast approximation algorithms for $\left(K, c_{\rm min} ,\alpha \right)$-coherent matrices. 

\subsection{Sublinear-Time Approximation Techniques}
\label{appsec:FastAlgs}

Consider the proof of Lemma~\ref{applem:Mest} with $c = 4$ for a given $m \times N$ $\left(K, c_{\rm min} ,\alpha \right)$-coherent matrix $\mathcal{M}$, $k \in \left[1,K \cdot c^2_{\rm min} / \alpha \right] \cap \mathbbm{N}$, and $\epsilon \in (0,1]$.  Let $\vec{x} \in \mathbbm{C}^N$ and suppose that 
\begin{equation}
x_n > 2 \delta = 2 \left( \frac{\epsilon \left\| \vec{x} - \vec{x}^{\rm opt}_{(k/\epsilon)} \right\|_1}{k} \right)
\label{appequ:DeltaDef}
\end{equation}
for some $n \in [0,N) \cap \mathbbm{N}$.  We will begin this section by quickly demonstrating a means of identifying $n$ using only the measurements $\mathcal{M} \vec{x} \in \mathbbm{C}^m$ together with some additional linear measurements based on a modification of our incoherent matrix $\mathcal{M}$.  This technique, first utilized in \cite{MuthuCS}, will ultimately allow us to develop the sublinear-time approximation schemes we seek.  However, we require several definitions before we can continue with our demonstration.

Let $\mathcal{A} \in \mathbbm{R}^{m \times N}$ and $\mathcal{C} \in \mathbbm{R}^{\tilde{m} \times N}$ be real matrices.  Then, their row tensor product, $\mathcal{A} \circledast \mathcal{C}$, is defined to be the $\left( m \cdot \tilde{m} \right) \times N$ matrix whose entries are given by
$$\left(\mathcal{A} \circledast \mathcal{C}\right)_{i,j} = \mathcal{A}_{i \textrm{~mod~m},j} \cdot \mathcal{C}_{\frac{i - i \textrm{~mod~m}}{m},j}.$$
In this section, we will use the row tensor product of $\mathcal{M}$ with the $\left( 1 + \lceil \log_2 N \rceil \right) \times N$ \textit{bit test matrix} \cite{MuthuCS,MyGroupTest} to help us identify $n$ from Equation~\ref{appequ:DeltaDef}.\footnote{We could also use the number theoretic $\mathcal{N}_{\lambda,s_1}$ matrices defined in Section 5 of \cite{ImpFourier} here in place of the bit test matrix.  More generally, any $1$-disjunct matrix with an associated fast decoding algorithm could replace the bit test matrix throughout this section.  Note that a fast $O(t)$-time binary tree decoder can be built for any $t \times N$ $1$-disjunct matrix anytime one has access to $\Omega(N)$ memory.}  The $\left( 1 + \lceil \log_2 N \rceil \right) \times N$ bit test matrix, $\mathcal{B}_N$, is defined by  
\begin{equation}
\left(\mathcal{B}_N\right)_{i,j} = \left\{ \begin{array}{ll} 1 & \textrm{if}~ i=0 \\ (i-1)^{\rm th} ~\textrm{bit in the binary expansion of}~j & \textrm{if}~i \in \left[1,\lceil \log_2 N \rceil \right]
\end{array} \right.
\label{appdef:BitTest}
\end{equation}
for $0 \leq i \leq \lceil \log_2 N \rceil$ and $ 0 \leq j < N$.
For example, $\mathcal{B}_8$ has the form
$$\mathcal{B}_8 = \left( \begin{array}{llllllll} 1 & 1 & 1 & 1 & 1 & 1 & 1 & 1 \\ 
0 & 1 & 0 & 1 & 0 & 1 & 0 & 1 \\ 
0 & 0 & 1 & 1 & 0 & 0 & 1 & 1 \\ 
0 & 0 & 0 & 0 & 1 & 1 & 1 & 1 \\ 
\end{array} \right).$$
We will now demonstrate that $\left(\mathcal{M} \circledast \mathcal{B}_N\right) \vec{x}$ contains enough information for us to identify any $n \in [0,N) \cap \mathbbm{N}$ satisfying Equation~\ref{appequ:DeltaDef}.

Notice that $\mathcal{M} \circledast \mathcal{B}_N$ contains $\mathcal{M}$ as a submatrix.  This is due to the first row of all ones in $\mathcal{B}_N$.  Similarly, the second row of $\mathcal{B}_N$ ensures that $\mathcal{M} \circledast \mathcal{B}_N$ will contain another $m \times N$ submatrix which is identical to $\mathcal{M}$, except with all of its even columns zeroed out.  We will refer to this $m \times N$ submatrix of $\mathcal{M} \circledast \mathcal{B}_N$ as $\mathcal{M}_{\rm odd}$.  We can see that
$$\mathcal{M}_{\rm odd} ~=~ \mathcal{M} \circledast \left(\mathcal{B}_N \right)_1 ~=~\left( \vec{0}~~\mathcal{M}_{\cdot,1}~~\vec{0}~~\mathcal{M}_{\cdot,3}~~\vec{0}~~\mathcal{M}_{\cdot,5}~~\vec{0}~~\dots \right).$$
Furthermore, we define 
$$\mathcal{M}_{\rm even} ~:=~  \mathcal{M} ~-~ \mathcal{M} \circledast \left(\mathcal{B}_N \right)_1 ~=~ \mathcal{M} ~-~ \mathcal{M}_{\rm odd}.$$  Clearly, if we are given $\left(\mathcal{M} \circledast \mathcal{B}_N\right) \vec{x}$, we will also have $\mathcal{M} \vec{x}, ~\mathcal{M}_{\rm odd} \vec{x}$, and $\mathcal{M}_{\rm even} \vec{x} \in \mathbbm{C}^m$.  We can use this information to determine whether $n$ from Equation~\ref{appequ:DeltaDef} is even or odd as follows.

Lemma~\ref{applem:Mest} with $c = 4$ guarantees that more than $K / 2$ distinct elements of $\mathcal{M} \vec{x} \in \mathbbm{C}^m$ will be of the form 
\begin{equation}
\left( \mathcal{M} \vec{x} \right)_j ~=~ x_n \cdot \mathcal{M}_{j,n} + (\mathcal{M'}(m,n))_j \cdot \vec{y}
\label{appequ:DecodeProp}
\end{equation}
for some $j \in [1,m]$ and $\vec{y} \in \mathbbm{C}^{N-1}$ with $\left| (\mathcal{M'}(m,n))_j \cdot \vec{y} \right| \leq c_{\rm min}\delta$ (see Equations~\ref{appdef:M'sub} and~\ref{appequ:DeltaDef} together with the proof of Lemma~\ref{applem:Mest}).  Suppose $n$ is odd.  Then, for each $j$ satisfying Equation~\ref{appequ:DecodeProp}, we will have
$$\left| \left( \mathcal{M}_{\rm even} \vec{x} \right)_j \right| ~=~ \left| \left( \mathcal{M'}_{\rm even}(m,n) \vec{y} \right)_j \right| ~\leq~c_{\rm min}\delta ~<~ |x_n| \cdot \mathcal{M}_{j,n} - \left| (\mathcal{M'}_{\rm odd}(m,n))_j \cdot \vec{y} \right| ~\leq~ \left| \left( \mathcal{M}_{\rm odd} \vec{x} \right)_j \right|.$$
Similarly, if $n$ is even, then for each such $j$ we will have
$$\left| \left( \mathcal{M}_{\rm odd} \vec{x} \right)_j \right| ~=~ \left| \left( \mathcal{M'}_{\rm odd}(m,n) \vec{y} \right)_j \right| ~\leq~c_{\rm min}\delta ~<~ |x_n| \cdot \mathcal{M}_{j,n} - \left| (\mathcal{M'}_{\rm even}(m,n))_j \cdot \vec{y} \right| ~\leq~ \left| \left( \mathcal{M}_{\rm even} \vec{x} \right)_j \right|.$$
Therefore, we can correctly determine $n ~{\rm mod}~ 2$ by comparing $\left| \left( \mathcal{M}_{\rm odd} \vec{x} \right)_j \right|$ with $\left| \left( \mathcal{M}_{\rm even} \vec{x} \right)_j \right|$ whenever both Equations~\ref{appequ:DeltaDef} and~\ref{appequ:DecodeProp} hold.  Of course, there is nothing particularly special about the lowest order bit of the binary representation of $n$.  More generally, we can correctly determine the $i^{\rm th}$ bit of $n \in [0,N) \cap \mathbbm{N}$ by comparing $\left| \left( \mathcal{M} \circledast \left(\mathcal{B}_N \right)_{i+1} \vec{x} \right)_j \right|$ with $\left| \left[ \left( \mathcal{M} ~-~ \mathcal{M} \circledast \left(\mathcal{B}_N \right)_{i+1} \right) \vec{x} \right]_j \right|$ whenever both Equations~\ref{appequ:DeltaDef} and~\ref{appequ:DecodeProp} hold.

We now know that we can correctly determine $n$ whenever both Equations~\ref{appequ:DeltaDef} and~\ref{appequ:DecodeProp} hold by finding its binary representation one bit at a time.  Furthermore, Lemma~\ref{applem:Mest} with $c \geq 4$ guarantees that more than $K / 2$ of the $j \in [1,m]$ will satisfy Equation~\ref{appequ:DecodeProp} for any given $n$.  Hence, we can correctly reconstruct every $n$ for which Equation~\ref{appequ:DeltaDef} holds more than $K / 2$ times by attempting to decode its binary representation for all $j \in [1,m]$.  Utilizing these methods together with ideas from Section~\ref{appsec:MProps}, we obtain Algorithm~\ref{appalg:reconstruct}.

\begin{algorithm}[tb]
\begin{algorithmic}[1]
\caption{$\proc{Approximate} ~\vec{x}$} \label{appalg:reconstruct}
\STATE \textbf{Input: An $m \times N$ $\left(K, c_{\rm min} ,\alpha \right)$-coherent matrix, $\mathcal{M}$, and $\left(\mathcal{M} \circledast \mathcal{B}_N\right) \vec{x} \in \mathbbm{C}^{m \lceil \log_2 N \rceil + m}$} 
\STATE \textbf{Output: $\vec{z}_S \approx \vec{x}_S$, an approximation to } $\vec{x}^{\rm~ opt}_k$
\STATE Initialize multiset $S \leftarrow \emptyset, ~\vec{z} \leftarrow \vec{0}_N, ~\vec{b} \leftarrow \vec{0}_{\lceil \log_2 N \rceil}$
\begin{center}
{\sc Identify All $n \in [0,N) \cap \mathbbm{N}$ that Satisfy Equation~\ref{appequ:DeltaDef}}
\end{center}
\FOR {$j$ from $1$ to $m$}
	\FOR {$i$ from $0$ to $\lceil \log_2 N \rceil-1$}
		\IF {$\left| \left( \mathcal{M} \circledast \left(\mathcal{B}_N \right)_{i+1} \vec{x} \right)_j \right| ~>~ \left| \left( \mathcal{M} \vec{x} ~-~ \mathcal{M} \circledast \left(\mathcal{B}_N \right)_{i+1} \vec{x} \right)_j \right|$}
			\STATE $b_i \leftarrow 1$
		\ELSE
			\STATE $b_i \leftarrow 0$
		\ENDIF
	\ENDFOR
	\STATE $n \leftarrow \sum^{\lceil \log_2 N \rceil-1}_{i=0} b_i 2^i$
	\STATE $S \leftarrow S \cup \{ n \}$
\ENDFOR \\
\begin{center}
{\sc Estimate $\vec{x}_S \approx \vec{x}^{\rm~ opt}_k$ Using Lemma~\ref{applem:Mest}}
\end{center}
\FOR {\textbf{each} $n$ value belonging to $S$ with multiplicity $> \frac{K}{2}$}
	\STATE $\mathbbm{Re}\left\{ z_n \right\} \leftarrow \textrm{median~of~multiset} \left\{ \mathbbm{Re} \left\{ \left(\mathcal{M}(K,n) \cdot \vec{x}_{} \right)_h / \left(\mathcal{M}(K,n)\right)_{h,n} \right\}~\big|~1 \leq h \leq K \right\}$
	\STATE $\mathbbm{Im}\left\{ z_n \right\} \leftarrow \textrm{median~of~multiset} \left\{ \mathbbm{Im} \left\{ \left(\mathcal{M}(K,n) \cdot \vec{x}_{} \right)_h / \left(\mathcal{M}(K,n)\right)_{h,n} \right\}~\big|~1 \leq h \leq K \right\}$
\ENDFOR
\STATE Sort nonzero $\vec{z}$ entries by magnitude so that $|z_{n_1}| \geq |z_{n_2}| \geq |z_{n_3}| \geq \dots$
\STATE $S \leftarrow \{ n_1, n_2, \dots, n_{2k} \}$
\STATE Output: $\vec{z}_{S}$
\end{algorithmic}
\end{algorithm}

In light of the preceding discussion, we can see that Algorithm~\ref{appalg:reconstruct} will be guaranteed to identify all $n \in [0,N) \cap \mathbbm{N}$ that satisfy Equation~\ref{appequ:DeltaDef} at least $\frac{K}{2}$ times each.  Lemma~\ref{applem:Mest} can then be used to estimate $x_n$ for each of these $n$ values as previously discussion in Section~\ref{appsec:MProps}.  The end result is that all relatively large entries in $\vec{x}$ will be identified and accurately estimated.  By formalizing the discussion above, we obtain the following result, the proof of which is analogous to the proof of Theorem 7 in Section 5 of \cite{ImpFourier}.

\begin{thm}
Suppose $\mathcal{M}$ is an $m \times N$ $\left(K, c_{\rm min} ,\alpha \right)$-coherent matrix.  Furthermore, let $\epsilon \in (0,1]$, $k \in \left[1,K \cdot \frac{\epsilon c^2_{\rm min}}{4 \alpha} \right) \cap \mathbbm{N}$, and $\vec{x} \in \mathbbm{C}^{N}$.  Then, Algorithm~\ref{appalg:reconstruct} will output a $\vec{z}_{S} \in \mathbbm{C}^N$ satisfying
$$ \left\| \vec{x} - \vec{z}_{S} \right\|_2 ~\leq~ \left\| \vec{x} - \vec{x}^{\rm~ opt}_k \right\|_2 + \frac{22 \epsilon \left\| \vec{x} - \vec{x}^{\rm opt}_{(k/\epsilon)} \right\|_1}{\sqrt{k}}.$$
Algorithm~\ref{appalg:reconstruct} can be implemented to run in $O\left( m \log N \right)$ time.
\label{appthm:RecoverAlg}
\end{thm}

The runtime of Algorithm~\ref{appalg:reconstruct} can be accounted for as follows:  Lines 4 through 14 can be implemented to run in $O(m \log N)$ time since their execution time will be dominated by the time required to read each entry of $\left(\mathcal{M} \circledast \mathcal{B}_N\right) \vec{x} \in \mathbbm{C}^{m \lceil \log_2 N \rceil + m}$.  Counting the multiplicity of the $O(m)$ entries in $S$ in line 15 can be done by sorting $S$ in $O(m \log m)$ time, followed by one $O(m)$-time scan of the sorted data.  Lines 16 and 17 will each be executed a total of $O(m / K)$ times apiece.  Furthermore, lines 16 and 17 can each be executed in $O(K \log K)$ time assuming that each $\mathcal{M}(K,n)$ submatrix is known in advance.\footnote{If each $\mathcal{M}(K,n)$ submatrix is not computed once in advance this runtime will be $O(m \log m)$ instead of $O(K \log K)$.}  Thus, the total runtime of lines 4 through 18 will also be $O(m \log N)$.  Finally, line 19 requires that at most $O(m/K)$ items be sorted, which can likewise be accomplished in $O(m \log N)$ time.  Therefore, the total runtime of Algorithm~\ref{appalg:reconstruct}  will be $O(m \log N)$.

\end{document}